\documentclass[11pt]{article}
\usepackage{amsmath,amssymb}
\usepackage{setspace,subfig}
\usepackage[hidelinks]{hyperref}
\usepackage{color}
\usepackage{tkz-graph}
\usepackage{pgf}
\usepackage{setspace}
\usepackage{tikz}
\usepackage{times}
\usetikzlibrary{arrows,shapes.arrows,shapes.geometric,
	shapes.multipart,backgrounds,decorations.pathmorphing,positioning,fit,automata}
	\usepackage{verbatim}

\usepackage{amsthm}
\usepackage{enumitem}
\usepackage{amssymb}
\usepackage{natbib}
\newcommand{\bm}[1]{\mbox{\boldmath{$#1$}}}

\usepackage{caption} 

\usepackage{graphicx}
\newcommand{\ind}{\rotatebox[origin=c]{90}{$\models$}}

\newtheoremstyle{note}
{8pt}
{8pt}
{}
{}
{\bfseries}
{:}
{.5em}
{}

\theoremstyle{note}

\newtheorem{remark}{Remark}

\newtheorem{definition}{Definition}

\newtheorem{prop}{Proposition}
\allowdisplaybreaks

\usepackage[left=1in,right=1in,top=1in,bottom=1in]{geometry}
\usepackage{algorithm}
\date{}
\usepackage{tikz}

\def\bSig\mathbf{\Sigma}

\usepackage[figuresright]{rotating}
\usepackage[utf8]{inputenc}
\usepackage{authblk}

\begin{document}

	\title{Conditional Independence Beyond Domain Separability \\
\Large{   Discussion of \cite{engelke2018graphical}}}
	\author[1]{Yuexia Zhang} 
	\author[1,2]{Linbo Wang\thanks{E-mail address for correspondence: linbo.wang@utoronto.ca}}
	\affil[1]{Department of Computer and Mathematical Sciences, University of Toronto Scarborough}
	\affil[2]{Department of Statistical Sciences, University of Toronto}
	\clearpage \maketitle

	\section{A general definition of conditional independence}

We congratulate Engelke and Hitz on a thought-provoking paper on graphical models for extremes. A key contribution of the paper is the introduction of a novel definition of conditional independence for a multivariate Pareto distribution. Here, we outline a proposal  for independence and conditional independence of general random variables whose support is a general set $\Omega$ in $\mathbb{R}^d$. Our proposal includes the authors' definition of conditional independence, and the analogous definition of independence as special cases. By making our proposal independent of the context of extreme value theory, we  highlight the importance of the authors' contribution beyond this particular context.

\begin{definition}
    Suppose that   $\mathbf Y=(\mathbf Y_A,\mathbf Y_B, \mathbf Y_C)$ is a random vector with support $\Omega$ in $\mathbb{R}^d$, where  $A, B,C$ are disjoint sets whose union is $\{1,\ldots,d\}$.  Let $U\times V$ denote the Cartesian product of $U$ and $V$.
    \begin{itemize}
        \item[(a)] We say $\mathbf Y_A$ is \emph{conditionally outer  independent} of $\mathbf Y_C$ given $\mathbf Y_B$ if there exists a random vector $\mathbf W = (\mathbf W_A, \mathbf W_B, \mathbf W_C)$ with  support $L_A \times L_B \times L_C$ in $\mathbb{R}^d$ such that (i) $\Omega \subset L_A \times L_B \times  L_C$; (ii) $(\mathbf W \mid \mathbf W\in \Omega) \,{\buildrel d \over =}\, \mathbf Y$; (iii) $\mathbf W_A \ind \mathbf W_C\mid \mathbf W_B$. In this case, we write $\mathbf Y_A \ind_o \mathbf Y_C\mid \mathbf Y_B.$
        
        If $B = \emptyset$, we say $\mathbf Y_A$ is \emph{outer  independent} of $\mathbf Y_C$, denoted as $\mathbf Y_A \ind_o \mathbf Y_C$.
        \item[(b)] We say $\mathbf Y_A$ is \emph{conditionally inner independent} of $\mathbf Y_C$ given $\mathbf Y_B$ if for any $S_A \times S_B \times S_C \subset \Omega$ such that $S_k$ is a measurable subset of $\mathbb{R}^{dim(Y_k)}, k\in\{A,B,C\}$ and $P(\mathbf Y\in S_A \times S_B \times S_C)>0$, we have $\mathbf Y_A \ind \mathbf Y_C\mid (\mathbf Y_B, \mathbf Y \in S_A\times S_B \times S_C)$. In this case, we write $\mathbf Y_A \ind_i \mathbf Y_C\mid \mathbf Y_B.$
       
       If $B = \emptyset$, we say $\mathbf Y_A$ is \emph{inner independent} of $\mathbf Y_C$, denoted as $\mathbf Y_A \ind_i \mathbf Y_C$.
    \end{itemize}
\end{definition}

\begin{prop}
\label{prop:ind}
    Suppose $\Omega = [0,\infty)^d \setminus [0,1]^{d}$ as in EH's case. Then 
   $$ \mathbf Y_A \ind_o \mathbf Y_C \mid \mathbf Y_B \quad \Leftrightarrow \quad \mathbf Y_A \ind_i \mathbf Y_C  \mid \mathbf Y_B \quad \Leftrightarrow \quad \mathbf Y_A \ind_e \mathbf Y_C \mid \mathbf Y_B.$$
   where $\ind_e$ denotes the notion of conditional independence introduced by \cite{engelke2018graphical}.
In particular, if $B = \emptyset$, then
   $$ \mathbf Y_A \ind_o \mathbf Y_C \quad \Leftrightarrow \quad \mathbf Y_A \ind_i \mathbf Y_C \quad \Leftrightarrow \quad \mathbf Y_A \ind_e \mathbf Y_C.$$
\end{prop}

\begin{remark}
 We do not place any distributional assumptions on $\mathbf Y$ in Proposition \ref{prop:ind}.
\end{remark}

 Engelke and Hitz showed that if $\mathbf Y$ is multivariate Pareto and admits a positive and continuous density, then $\mathbf Y_A \not\!\perp\!\!\!\perp_e \mathbf Y_C.$ This does not rule out the possibility of $\mathbf Y_A \ind_e \mathbf Y_C$ for general Pareto distributions.  For example, consider two independent standard Pareto distributions $X_1$ and $X_2$.
Following eqn. (6) in \cite{engelke2018graphical}, all the probability mass of $\mathbf Y$ lies on $(1,\infty)\times \{0\}$ and $\{0\}\times (1,\infty)$ so that it does not admit a density with respect to Lebesgue measure. Nevertheless  $Y_1 \ind_e Y_2$.
This observation is generalized in Proposition \ref{prop:eh}.

\begin{prop}
	\label{prop:eh}
Suppose $\Omega = \left[ ( [0,\infty)^{|A|} \setminus [0,1]^{|A|}) \times \{0\}^{|B|}  \times \{0\}^{|C|}\right] \cup \left[\{0\}^{|A|}  \times ([0,\infty)^{|B|} \setminus [0,1]^{|B|}) \times  \{0\}^{|C|} \right] \\ \cup \left[\{0\}^{|A|}  \times \{0\}^{|B|}  \times ([0,\infty)^{|C|} \setminus [0,1]^{|C|} ) \right],$ where $|A|$ is the cardinality of set $A$.  Then 
	$$\mathbf Y_A \ind_o\mathbf Y_C \mid \mathbf Y_B \quad \Leftrightarrow \quad \mathbf Y_A \ind_i \mathbf Y_C  \mid \mathbf Y_B \quad \Leftrightarrow \quad \mathbf Y_A \ind_e \mathbf Y_C \mid \mathbf Y_B.$$
	In particular, if $B = \emptyset$, then
	$$\mathbf Y_A \ind_o \mathbf Y_C \quad \Leftrightarrow \quad \mathbf Y_A \ind_i \mathbf Y_C \quad \Leftrightarrow \quad \mathbf Y_A \ind_e \mathbf Y_C.$$
\end{prop}

\section{Proofs of propositions}

\subsection{Proof of Proposition \ref{prop:ind} in the case where  $B=\emptyset$}

\begin{proof}
We need to show
\begin{equation}
    \label{eqn:claim1}
     \mathbf Y_A \ind_o \mathbf Y_C \Rightarrow \mathbf Y_A \ind_i \mathbf Y_C \Rightarrow \mathbf Y_A \ind_e \mathbf Y_C \Rightarrow \mathbf Y_A \ind_o \mathbf Y_C.  
\end{equation}
  
The first two claims are straightforward. We now prove the third claim in \eqref{eqn:claim1}. 

Let $H_k^1\subset L_k^1=[0,1]^{|k|}$, $H_k^2\subset L_k^2=[0,\infty)^{|k|}\setminus [0,1]^{|k|}$, $H_k=H_k^1\cup H_k^2$, $L_k=L_k^1\cup L_k^2$, $k\in \{A,C\}$. Obviously, $H_k^1\cap H_k^2=\emptyset$ and $L_k^1\cap L_k^2=\emptyset$. Besides, $\Omega = [0,\infty)^d \setminus [0,1]^{d} \subset L_A \times L_C$.

We let 
\begin{equation}
    \label{eqn:definition1}
        {P(\mathbf W_A \in H_A^1, \mathbf W_C\in H_C^1)} = \dfrac{P(\mathbf W_A \in H_A^1, \mathbf W_C\in L_C^2){P(\mathbf W_A\in L_A^2, \mathbf W_C\in H_C^1)}}{P(\mathbf W_A\in L_A^2, \mathbf W_C\in L_C^2)},
\end{equation}
and 
\begin{equation*}
P(\mathbf W_A\in H_A^i, \mathbf W_C\in H_C^j) = \lambda P (\mathbf Y_A \in H_A^i, \mathbf Y_C \in H_C^j), \quad (i,j)\neq (1,1),
\end{equation*}
where $\lambda$ is a normalizing constant.
Simple calculation yields
$$
        \lambda =\dfrac{1}{1+ \dfrac{P(\mathbf Y_A\in L_A^1, \mathbf Y_C \in L_C^2) P(\mathbf Y_A\in L_A^2, \mathbf Y_C\in L_C^1)}{P(\mathbf Y_A\in L_A^2,\mathbf Y_C\in L_C^2)}}.
$$

It is easy to see that $(\mathbf W \mid \mathbf W \in \Omega ) \,{\buildrel d \over =}\, \mathbf Y.$ 
We then prove
\begin{equation*}
    \mathbf W_A \ind \mathbf W_C.
\end{equation*}

If $k \in A$ and $\{\mathbf W_A \in H_A\}\cap \{\mathbf W_A: W_k>1\}=\emptyset$, then
\begin{equation}
\label{eqn:1}
P(\mathbf W_A\in H_A, \mathbf W_C\in H_C\mid W_k>1)=P(\mathbf W_A\in H_A \mid W_k>1)P(\mathbf W_C\in H_C \mid W_k>1)=0.
\end{equation}
If $k \in A$ and $\{\mathbf W_A \in H_A\}\cap \{\mathbf W_A: W_k>1\}\neq \emptyset$, then $\{\mathbf W_A \in H_A\}\cap \{\mathbf W_A: W_k>1\}=\{\mathbf W_A \in H_A^2, W_k>1\}$ and

\begin{flalign*}
&P(\mathbf W_A \in H_A,\mathbf W_C \in H_C, W_k>1) \\
=& \lambda P (\mathbf Y_A \in H_A,\mathbf Y_C \in H_C, Y_k>1) \\
=& \lambda  P (\mathbf Y_C\in H_C \mid  \mathbf Y_A\in H_A, Y_k > 1) P(\mathbf Y_A\in H_A, Y_k > 1) \\
=& P (\mathbf Y_C \in H_C \mid Y_k > 1) P(\mathbf W_A \in H_A, W_k>1) \quad (\mathbf Y_A \ind \mathbf Y_C \mid Y_k > 1).
\end{flalign*}
So we have  
\begin{equation}
\begin{aligned}
 \label{eqn:2}
    P(\mathbf W_C\in H_C \mid \mathbf W_A\in H_A, W_k>1) &= P (\mathbf Y_C \in H_C \mid Y_k > 1)\\
    &= \dfrac{P (\mathbf Y_C \in H_C , Y_k > 1)}{P(Y_k > 1)} \\
    &= \dfrac{P (\mathbf W_C \in H_C , W_k > 1)/\lambda}{P(W_k > 1)/\lambda} \\
    &= \dfrac{P (\mathbf W_C \in H_C , W_k > 1)}{P(W_k > 1)} \\
    &= P(\mathbf W_C \in H_C \mid W_k>1). 
\end{aligned}
\end{equation}
Thus, from \eqref{eqn:1} and \eqref{eqn:2}, we have for any $k \in A$, $\mathbf W_A \ind \mathbf W_C \mid W_k>1$. Similarly, we have for any $k \in C$, $\mathbf W_A \ind \mathbf W_C \mid W_k>1$. In summary, 
\begin{equation}
\label{eqn:ind1}
   \forall k \in \{1,\ldots,d\}: \mathbf W_A \ind \mathbf W_C \mid W_k>1. 
\end{equation}

We then show that  for any $H_A^1, H_C^1, H_A^2, H_C^2$, we have
\begin{equation}
\label{eqn:cons1}
\dfrac{P(\mathbf W_A \in H_A^1, \mathbf W_C\in H_C^1)}{P(\mathbf W_A\in H_A^2, \mathbf W_C\in H_C^1)} = \dfrac{P(\mathbf W_A \in H_A^1,\mathbf  W_C\in H_C^2)}{P(\mathbf W_A\in H_A^2, \mathbf W_C\in H_C^2)}.
\end{equation}

To prove \eqref{eqn:cons1}, note that 
\begin{flalign*}
& \ \ \dfrac{P(\mathbf W_A \in H_A^1, \mathbf W_C\in H_C^2)}{P(\mathbf W_A\in H_A^2, \mathbf W_C\in H_C^2)}P(\mathbf W_A\in H_A^2, \mathbf W_C\in H_C^1)\\ 
&= \dfrac{P(\mathbf W_A \in H_A^1 \mid \mathbf W_C\in H_C^2)}{P(\mathbf W_A\in H_A^2 \mid \mathbf W_C\in H_C^2)}P(\mathbf W_A\in H_A^2, \mathbf W_C\in H_C^1) \\
&= \dfrac{P(\mathbf W_A \in H_A^1 \mid \mathbf W_C\in H_C^2, \mathbf W_C\in L_C^2)}{P(\mathbf W_A\in H_A^2 \mid \mathbf W_C\in H_C^2, \mathbf W_C\in L_C^2)}P(\mathbf W_A\in H_A^2, \mathbf W_C\in H_C^1) \\
&= \dfrac{P(\mathbf W_A \in H_A^1 \mid \mathbf W_C\in L_C^2)}{P(\mathbf W_A\in H_A^2 \mid \mathbf W_C\in L_C^2)} P(\mathbf W_A\in H_A^2, \mathbf W_C\in H_C^1) \quad \text{due to } \eqref{eqn:ind1} \\
&= \dfrac{P(\mathbf W_A\in H_A^2, \mathbf W_C\in H_C^1)}{P(\mathbf W_A\in H_A^2, \mathbf W_C\in L_C^2)}  P(\mathbf W_A \in H_A^1,\mathbf  W_C\in L_C^2) \\
&= \dfrac{P(\mathbf W_A\in H_A^2, \mathbf W_A\in L_A^2, \mathbf W_C\in H_C^1)}{P(\mathbf W_A\in H_A^2, \mathbf W_A\in L_A^2, \mathbf W_C\in L_C^2)}  P(\mathbf W_A \in H_A^1, \mathbf W_C\in L_C^2) \\
&= \dfrac{P(\mathbf W_A \in L_A^2, \mathbf W_C\in H_C^1)}{P(\mathbf W_A \in L_A^2, \mathbf W_C \in L_C^2)}  P(\mathbf W_A \in H_A^1, \mathbf W_C \in L_C^2) \quad \text{due to } \eqref{eqn:ind1}\\
&=P(\mathbf W_A \in H_A^1, \mathbf W_C\in H_C^1)\quad \text{due to } \eqref{eqn:definition1}.
\end{flalign*}

We then show that 
\begin{equation}
    \label{eqn:3}
    P(\mathbf W_C \in L_C^2 \mid \mathbf W_A \in L_A^2) = P(\mathbf W_C \in L_C^2 \mid \mathbf W_A \in L_A^1).
\end{equation}
To prove \eqref{eqn:3}, note that due to \eqref{eqn:cons1},
$$
        \dfrac{P(\mathbf W_A \in L_A^1, \mathbf W_C \in L_C^1)}{P(\mathbf W_A \in L_A^2, \mathbf W_C \in L_C^1)}=\dfrac{P(\mathbf W_A \in L_A^1, \mathbf W_C \in L_C^2)}{P(\mathbf W_A \in L_A^2, \mathbf W_C \in L_C^2)}.
$$
Then
\begin{flalign*}
\text{LHS of } \eqref{eqn:3} &= \dfrac{P(\mathbf W_A \in L_A^2,\mathbf  W_C \in L_C^2)}{P(\mathbf W_A \in L_A^2)} \\
&= \dfrac{P(\mathbf W_A \in L_A^2, \mathbf W_C \in L_C^2)}{P(\mathbf W_A \in L_A^2,\mathbf W_C \in L_C^2) + P(\mathbf W_A \in L_A^2, \mathbf W_C \in L_C^1)}  \\
&= \dfrac{P(\mathbf W_A \in L_A^1, \mathbf W_C \in L_C^2)}{P(\mathbf W_A \in L_A^1,\mathbf W_C \in L_C^2) + P(\mathbf W_A \in L_A^1, \mathbf W_C \in L_C^1)}  \\
&= \dfrac{P(\mathbf W_A \in L_A^1, \mathbf W_C \in L_C^2)}{P(\mathbf W_A \in L_A^1)}  \\
&= P(\mathbf W_C \in L_C^2\mid  \mathbf W_A \in L_A^1) = \text{RHS of } \eqref{eqn:3}.
\end{flalign*}
Similarly we have 
\begin{flalign*}
P(\mathbf W_C \in L_C^1\mid \mathbf W_A \in L_A^2) &= P(\mathbf W_C \in L_C^1\mid \mathbf W_A \in L_A^1).
\end{flalign*}

We then show that 
\begin{equation*}
    P(\mathbf W_C\in H_C^2\mid \mathbf W_A \in L_A^2) = P(\mathbf W_C \in H_C^2 \mid\mathbf  W_A \in L_A^1).
\end{equation*}
To see this, due to \eqref{eqn:cons1}, 
\begin{flalign*}
\dfrac{P(\mathbf W_A \in L_A^1, \mathbf W_C\in H_C^1)}{P(\mathbf W_A \in L_A^2, \mathbf W_C\in H_C^1)} = \dfrac{P(\mathbf W_A \in L_A^1, \mathbf W_C \in L_C^2)}{P(\mathbf W_A \in L_A^2, \mathbf W_C \in L_C^2)}.
\end{flalign*}
Furthermore, due to \eqref{eqn:3},
$$
    \dfrac{P(\mathbf W_C\in H_C^1\mid \mathbf W_A \in L_A^1)}{P(\mathbf W_C\in H_C^1\mid \mathbf W_A \in L_A^2)} = \dfrac{P(\mathbf W_C \in L_C^2\mid \mathbf  W_A \in L_A^1)}{P(\mathbf W_C \in L_C^2\mid \mathbf W_A \in L_A^2)} = 1.
$$
Due to \eqref{eqn:ind1} and \eqref{eqn:cons1},
it follows directly that 
\begin{multline}
\label{eqn:4}
P(\mathbf W_C\in H_C^2\mid \mathbf W_A \in L_A^1) 
= P(\mathbf W_C\in H_C^2\mid \mathbf W_A \in L_A^2) \\
= P(\mathbf W_C\in H_C^2\mid \mathbf W_A \in H_A^2) 
= P(\mathbf W_C\in H_C^2).
\end{multline}
Similarly we have
\begin{equation}
\label{eqn:5}
P(\mathbf W_C\in H_C^1\mid \mathbf W_A \in L_A^1)= P(\mathbf W_C\in H_C^1 \mid \mathbf W_A \in H_A^2) = P(\mathbf W_C\in H_C^1).
\end{equation}

We still need to show 
\begin{flalign}
P(\mathbf W_C \in H_C^2 \mid \mathbf W_A \in H_A^1) &= P(\mathbf W_C \in H_C^2 \mid \mathbf W_A \in L_A^1), \label{eqn:6}\\
P(\mathbf W_C \in H_C^1 \mid \mathbf W_A \in H_A^1) &= P(\mathbf W_C \in H_C^1 \mid \mathbf W_A \in L_A^1).\label{eqn:7}
\end{flalign}

To show \eqref{eqn:6}, note that due to \eqref{eqn:ind1} and \eqref{eqn:cons1},

\begin{flalign*}
    & \ \ P(\mathbf W_C \in H_C^2\mid \mathbf W_A \in H_A^1)    \\
    &= \dfrac{P(\mathbf W_C \in H_C^2, \mathbf W_A \in H_A^1)}{P(\mathbf W_A \in H_A^1)} \\
    &= \dfrac{ \lambda  P (\mathbf Y_C \in H_C^2,\mathbf  Y_C\in L_C^2, \mathbf Y_A \in H_A^1)}{P(\mathbf W_A \in H_A^1, \mathbf W_C\in L_C^1) + P(\mathbf W_A \in H_A^1, \mathbf W_C\in L_C^2)} \\
    &= \dfrac{ \lambda   P(\mathbf Y_C \in H_C^2 \mid  \mathbf Y_C\in L_C^2) P(\mathbf Y_A \in H_A^1, \mathbf Y_C\in L_C^2)}{ \dfrac{P(\mathbf W_A \in H_A^1,  \mathbf W_C\in L_C^2){P(\mathbf W_A \in L_A^2, \mathbf W_C\in L_C^1)}}{P(\mathbf W_A \in L_A^2, \mathbf W_C\in L_C^2)} + P(\mathbf W_A \in H_A^1, \mathbf W_C\in L_C^2)} \\
    &= \dfrac{   P(\mathbf Y_C \in H_C^2 \mid \mathbf Y_C\in L_C^2) P(\mathbf W_A \in H_A^1, \mathbf W_C\in L_C^2)}{ \dfrac{P( \mathbf W_A \in H_A^1,  \mathbf W_C\in L_C^2){P(\mathbf W_A \in L_A^2,\mathbf  W_C\in L_C^1)}}{P(\mathbf W_A \in L_A^2, \mathbf W_C\in L_C^2)} + P(\mathbf W_A \in H_A^1, \mathbf W_C\in L_C^2)} \\
      &= \dfrac{   P(\mathbf Y_C \in H_C^2 \mid \mathbf  Y_C\in L_C^2) }{ \dfrac{{P(\mathbf W_A \in L_A^2, \mathbf W_C\in L_C^1)}}{P(\mathbf W_A \in L_A^2, \mathbf W_C\in L_C^2)} + 1} 
\end{flalign*}
does not depend on $H_A^1$.

To show \eqref{eqn:7}, note that due to \eqref{eqn:definition1} and \eqref{eqn:cons1},
\begin{flalign*}
& \ \ P(\mathbf W_C \in H_C^1 \mid \mathbf W_A \in H_A^1) \\ &= \dfrac{P(\mathbf W_A \in H_A^1,  \mathbf W_C \in H_C^1)}{P(\mathbf W_A \in H_A^1,  \mathbf W_C \in L_C^1) + P(\mathbf W_A \in H_A^1,\mathbf  W_C \in L_C^2)} \\
&= \dfrac{ \dfrac{P(\mathbf W_A \in H_A^1, \mathbf W_C \in L_C^2){P(\mathbf W_A \in L_A^2, \mathbf W_C \in H_C^1)}}{P(\mathbf W_A \in L_A^2,  \mathbf W_C \in L_C^2)} }{\dfrac{P(\mathbf W_A \in H_A^1, \mathbf W_C \in L_C^2){P(\mathbf W_A \in L_A^2, \mathbf  W_C \in L_C^1)}}{P(\mathbf W_A \in L_A^2,  \mathbf W_C \in L_C^2)} + P(\mathbf W_A \in H_A^1, \mathbf W_C \in L_C^2)} \\
&= \dfrac{ \dfrac{{P(\mathbf W_A \in L_A^2,  \mathbf W_C \in H_C^1)}}{P(\mathbf W_A \in L_A^2, \mathbf W_C \in L_C^2)} }{\dfrac{{P(\mathbf W_A \in L_A^2, \mathbf W_C \in L_C^1)}}{P(\mathbf W_A \in L_A^2, \mathbf W_C \in L_C^2)} + 1} 
\end{flalign*}
does not depend on $H_A^1$.

From \eqref{eqn:4}, \eqref{eqn:5}, \eqref{eqn:6} and \eqref{eqn:7}, we have
\begin{equation*}
        P(\mathbf W_C \in H_C \mid \mathbf W_A \in H_A)=P(\mathbf W_C \in H_C).
\end{equation*}
Hence $\mathbf W_A \ind \mathbf W_C$.

\end{proof}

\subsection{Proof of Proposition \ref{prop:ind} for a general $B$}
\begin{proof}

We need to show
\begin{equation}
    \label{eqn:claim2}
      \mathbf Y_A \ind_o \mathbf Y_C \mid \mathbf Y_B \Rightarrow\mathbf  Y_A \ind_i \mathbf Y_C \mid \mathbf Y_B \Rightarrow \mathbf Y_A \ind_e \mathbf Y_C \mid \mathbf Y_B \Rightarrow \mathbf Y_A \ind_o \mathbf Y_C \mid\mathbf  Y_B.  
\end{equation}
  
The first two claims are straightforward. We now prove the third claim in \eqref{eqn:claim2}.  

Let $H_k^1\subset L_k^1=[0,1]^{|k|}$, $H_k^2\subset L_k^2=[0,\infty)^{|k|}\setminus [0,1]^{|k|}$, $H_k=H_k^1\cup H_k^2$, $L_k=L_k^1\cup L_k^2$, $k\in \{A,B,C\}$. Obviously, $H_k^1\cap H_k^2=\emptyset$ and $L_k^1\cap L_k^2=\emptyset$. Besides, $\Omega = [0,\infty)^d \setminus [0,1]^{d} \subset L_A \times L_B \times L_C$.

We let 
\begin{multline}
\label{eqn:definition3}
        {P(\mathbf W_A \in H_A^1,\mathbf  W_B\in H_B^1, \mathbf W_C\in H_C^1)} \\
        = \dfrac{P(\mathbf W_A \in H_A^1, \mathbf W_B\in H_B^1, \mathbf W_C\in L_C^2){P(\mathbf W_A\in L_A^2, \mathbf W_B\in H_B^1, \mathbf W_C\in H_C^1)}}{P(\mathbf W_A\in L_A^2, \mathbf W_B\in H_B^1, \mathbf W_C\in L_C^2)},
\end{multline}
and for $(i, j, k)\neq (1, 1, 1)$,
\begin{equation*}
P(\mathbf W_A\in H_A^i, \mathbf W_B\in H_B^j, \mathbf W_C\in H_C^k) = \lambda P (\mathbf Y_A \in H_A^i, \mathbf Y_B\in H_B^j, \mathbf Y_C \in H_C^k),  
\end{equation*}
where $\lambda$ is a normalizing constant.
Simple calculation yields
$$
        \lambda =\dfrac{1}{1+ \dfrac{P(\mathbf Y_A\in L_A^1,\mathbf  Y_B\in L_B^1, \mathbf Y_C \in L_C^2) P(\mathbf Y_A\in L_A^2, \mathbf Y_B\in L_B^1, \mathbf Y_C\in L_C^1)}{P(\mathbf Y_A\in L_A^2, \mathbf Y_B\in L_B^1, \mathbf Y_C\in L_C^2)}}.
$$

It is easy to see that $(\mathbf  W \mid  \mathbf  W \in \Omega ) \,{\buildrel d \over =}\, \mathbf  Y.$ 
We then prove
\begin{equation*}
   \mathbf  W_A \ind\mathbf  W_C \mid \mathbf W_B.
\end{equation*}

If $k \in A$ and $\{\mathbf W_A \in H_A\}\cap \{\mathbf W_A: W_k>1\}=\emptyset$, then
\begin{multline}
\label{eqn:37}
P(\mathbf W_A\in H_A, \mathbf W_C\in H_C\mid \mathbf W_B\in H_B, W_k>1)=P(\mathbf W_A\in H_A \mid \mathbf W_B\in H_B, W_k>1)\\
\times P(\mathbf W_C\in H_C \mid \mathbf W_B\in H_B, W_k>1)=0.
\end{multline}
If $k \in A$ and $\{\mathbf W_A \in H_A\}\cap \{\mathbf W_A: W_k>1\}\neq \emptyset$, then $\{\mathbf W_A \in H_A\}\cap \{\mathbf W_A: W_k>1\}=\{\mathbf W_A \in H_A^2, W_k>1\}$.
Note that
\begin{flalign*}
&P(\mathbf W_A \in H_A, \mathbf W_B \in H_B, \mathbf W_C \in H_C, W_k>1) \\
&= \lambda P (\mathbf Y_A \in H_A, \mathbf Y_B \in H_B, \mathbf Y_C \in H_C, Y_k>1) \\
&= \lambda  P (\mathbf Y_C\in H_C \mid \mathbf  Y_A\in H_A, \mathbf Y_B \in H_B, Y_k > 1) P( \mathbf Y_A\in H_A, \mathbf Y_B \in H_B, Y_k > 1) \\
&= P (\mathbf Y_C \in H_C \mid\mathbf  Y_B \in H_B, Y_k > 1) P(\mathbf W_A \in H_A, \mathbf W_B \in H_B, W_k>1) \quad (\mathbf Y_A \ind \mathbf Y_C \mid\mathbf  Y_B, Y_k > 1).
\end{flalign*}
So we have  
\begin{equation}
\begin{aligned}
 \label{eqn:38}
    P(\mathbf W_C\in H_C \mid \mathbf W_A\in H_A, \mathbf W_B \in H_B, W_k>1) &= P (\mathbf Y_C \in H_C \mid\mathbf  Y_B \in H_B, Y_k > 1)\\
    &= \dfrac{P (\mathbf Y_C \in H_C, \mathbf Y_B \in H_B, Y_k > 1)}{P(\mathbf Y_B \in H_B, Y_k > 1)} \\
    &= \dfrac{P (\mathbf W_C \in H_C, \mathbf W_B \in H_B, W_k > 1)/\lambda}{P(\mathbf W_B \in H_B, W_k > 1)/\lambda} \\
    &= \dfrac{P (\mathbf W_C \in H_C, \mathbf W_B \in H_B, W_k > 1)}{P(\mathbf W_B \in H_B, W_k > 1)} \\
    &= P(\mathbf W_C \in H_C \mid \mathbf W_B \in H_B, W_k>1). 
\end{aligned}
\end{equation}
Thus, from \eqref{eqn:37} and \eqref{eqn:38}, we have for any $k \in A$, $\mathbf W_A \ind \mathbf W_C \mid \mathbf W_B, W_k>1$.
Similarly, we have for any $k \in C$, $\mathbf W_A \ind \mathbf W_C \mid \mathbf W_B, W_k>1$.

If $k \in B$ and $\{\mathbf W_B \in H_B\}\cap \{\mathbf W_B: W_k>1\}=\emptyset$, then according to the definition in \citet[\S 7.4]{shorack2017probability}, 
\begin{equation}
    \label{eqn:39}
    P(\mathbf W_C \in H_C \mid \mathbf W_A \in H_A, \mathbf W_B \in H_B,  W_k>1)=P(\mathbf W_C \in H_C).
\end{equation}
If $k \in B$ and $\{\mathbf W_B \in H_B\}\cap \{\mathbf W_B: W_k>1\}\neq \emptyset$, then $\{\mathbf W_B\in H_B\}\cap \{\mathbf W_B: W_k>1\}=\{\mathbf W_B \in H_B^2, W_k>1\}$. Besides, similar to the proof of  \eqref{eqn:38}, we have 
\begin{equation}
    \label{eqn:40}
    P(\mathbf W_C \in H_C \mid \mathbf W_A \in H_A, \mathbf W_B \in H_B,  W_k>1)=P(\mathbf W_C \in H_C \mid \mathbf W_B \in H_B,  W_k>1 ).
\end{equation}
Thus, from \eqref{eqn:39} and \eqref{eqn:40}, we have for any $k \in B$, $\mathbf W_A \ind \mathbf W_C \mid \mathbf W_B, W_k>1$.

In summary,
\begin{equation}
    \label{eqn:ind3}
   \forall k \in \{1,\ldots,d\}: \mathbf W_A\ind \mathbf W_C \mid \mathbf W_B, W_k>1.
\end{equation}

From  \eqref{eqn:ind3}, obviously, we have
\begin{equation}
   \label{eqn:41} 
    P(\mathbf W_C \in H_C \mid \mathbf W_A \in H_A,\mathbf  W_B \in H_B^2)=P(\mathbf W_C \in H_C \mid \mathbf  W_B \in H_B^2).
\end{equation}

We then show that  for any $H_A^1, H_B^1, H_C^1, H_A^2, H_C^2$, we have
\begin{equation}
\label{eqn:cons3}
\dfrac{P(\mathbf W_A \in H_A^1, \mathbf W_C\in H_C^1\mid \mathbf W_B\in H_B^1 )}{P(\mathbf W_A\in H_A^2, \mathbf W_C\in H_C^1\mid \mathbf W_B \in H_B^1)} = \dfrac{P(\mathbf W_A \in H_A^1, \mathbf W_C\in H_C^2\mid \mathbf W_B \in H_B^1)}{P(\mathbf W_A\in H_A^2,\mathbf  W_C\in H_C^2\mid \mathbf W_B \in H_B^1)}.
\end{equation}

To prove \eqref{eqn:cons3}, note that 
\begin{flalign*}
& \ \ \dfrac{P(\mathbf W_A \in H_A^1, \mathbf W_C\in H_C^2\mid \mathbf W_B \in H_B^1)}{P(\mathbf W_A\in H_A^2,\mathbf  W_C\in H_C^2\mid \mathbf W_B \in H_B^1)}P(\mathbf W_A\in H_A^2, \mathbf W_C\in H_C^1\mid \mathbf W_B\in H_B^1)\\ &= \dfrac{P(\mathbf W_A \in H_A^1 \mid \mathbf W_C\in H_C^2,\mathbf W_B\in H_B^1)}{P(\mathbf W_A\in H_A^2 \mid \mathbf W_C\in H_C^2,\mathbf W_B\in H_B^1)}P(\mathbf W_A\in H_A^2,\mathbf  W_C\in H_C^1\mid \mathbf W_B\in H_B^1) \\
&= \dfrac{P(\mathbf W_A \in H_A^1 \mid \mathbf W_C\in L_C^2, \mathbf W_B\in H_B^1)}{P(\mathbf W_A\in H_A^2 \mid \mathbf W_C\in L_C^2, \mathbf W_B\in H_B^1)} P(\mathbf W_A\in H_A^2, \mathbf W_C\in H_C^1\mid \mathbf W_B\in H_B^1) \quad \text{due to } \eqref{eqn:ind3} \\
&= \dfrac{P(\mathbf W_A\in H_A^2, \mathbf W_C\in H_C^1 \mid \mathbf W_B\in H_B^1)}{P(\mathbf W_A\in H_A^2, \mathbf W_C\in L_C^2 \mid \mathbf W_B\in H_B^1)}  P(\mathbf W_A \in H_A^1, \mathbf W_C\in L_C^2 \mid \mathbf W_B\in H_B^1) \\
&= \dfrac{P(\mathbf W_A \in L_A^2,\mathbf  W_C\in H_C^1\mid \mathbf W_B\in H_B^1)}{P(\mathbf W_A \in L_A^2, \mathbf W_C \in L_C^2 \mid \mathbf W_B\in H_B^1)}  P(\mathbf W_A \in H_A^1,\mathbf  W_C \in L_C^2\mid \mathbf W_B\in H_B^1) \quad \text{due to } \eqref{eqn:ind3}\\
&=P(\mathbf W_A \in H_A^1, \mathbf W_C\in H_C^1\mid \mathbf W_B\in H_B^1 )\quad \text{due to } \eqref{eqn:definition3}.
\end{flalign*}

We then show that 
\begin{equation}
    \label{eqn:42}
    P(\mathbf W_C \in L_C^2 \mid \mathbf W_A \in L_A^2, \mathbf W_B\in H_B^1) = P(\mathbf W_C \in L_C^2 \mid \mathbf W_A \in L_A^1, \mathbf W_B\in H_B^1).
\end{equation}
To prove \eqref{eqn:42}, note that due to \eqref{eqn:cons3},
$$
        \dfrac{P(\mathbf W_A \in L_A^1, \mathbf W_C \in L_C^1\mid \mathbf W_B\in H_B^1)}{P(\mathbf W_A \in L_A^2, \mathbf W_C \in L_C^1\mid \mathbf W_B\in H_B^1)}=\dfrac{P(\mathbf W_A \in L_A^1, \mathbf W_C \in L_C^2\mid \mathbf W_B\in H_B^1)}{P(\mathbf W_A \in L_A^2, \mathbf W_C \in L_C^2\mid \mathbf W_B\in H_B^1)}.
$$
Then
\begin{flalign*}
\text{LHS of } \eqref{eqn:42} &= \dfrac{P(\mathbf W_A \in L_A^2, \mathbf W_C \in L_C^2, \mathbf W_B\in H_B^1)}{P(\mathbf W_A \in L_A^2,\mathbf W_B\in H_B^1)} \\
&= \dfrac{P(\mathbf W_A \in L_A^2, \mathbf W_C \in L_C^2,\mathbf W_B\in H_B^1)}{P(\mathbf W_A \in L_A^2,\mathbf W_C \in L_C^2,\mathbf W_B\in H_B^1) + P(\mathbf W_A \in L_A^2, \mathbf W_C \in L_C^1, \mathbf W_B\in H_B^1)}  \\
&= \dfrac{P(\mathbf W_A \in L_A^1, \mathbf W_C \in L_C^2\mid \mathbf W_B\in H_B^1)}{P(\mathbf W_A \in L_A^1,\mathbf W_C \in L_C^2\mid \mathbf W_B\in H_B^1) + P(\mathbf W_A \in L_A^1, \mathbf W_C \in L_C^1 \mid \mathbf W_B\in H_B^1)}  \\
&= \dfrac{P(\mathbf W_A \in L_A^1, \mathbf W_C \in L_C^2\mid \mathbf W_B\in H_B^1)}{P(\mathbf W_A \in L_A^1\mid \mathbf W_B\in H_B^1)}  \\
&= P(\mathbf W_C \in L_C^2\mid  \mathbf W_A \in L_A^1, \mathbf W_B\in H_B^1) = \text{RHS of } \eqref{eqn:42}.
\end{flalign*}
Similarly we have 
\begin{flalign*}
P(\mathbf W_C \in L_C^1\mid \mathbf W_A \in L_A^2,\mathbf W_B\in H_B^1) &= P(\mathbf W_C \in L_C^1\mid \mathbf W_A \in L_A^1,\mathbf W_B\in H_B^1).
\end{flalign*}

We then show that 
\begin{equation*}
    P(\mathbf W_C\in H_C^2\mid \mathbf W_A \in L_A^2,\mathbf W_B\in H_B^1) = P(\mathbf W_C \in H_C^2 \mid \mathbf W_A \in L_A^1,\mathbf W_B\in H_B^1).
\end{equation*}
To see this, due to \eqref{eqn:cons3}, 
\begin{flalign*}
\dfrac{P(\mathbf W_A \in L_A^1, \mathbf W_C\in H_C^1\mid \mathbf W_B\in H_B^1)}{P(\mathbf W_A \in L_A^2, \mathbf W_C\in H_C^1\mid \mathbf W_B\in H_B^1)} = \dfrac{P(\mathbf W_A \in L_A^1, \mathbf W_C \in L_C^2\mid \mathbf W_B\in H_B^1)}{P(\mathbf W_A \in L_A^2, \mathbf W_C \in L_C^2\mid \mathbf W_B\in H_B^1)}.
\end{flalign*}
Furthermore, due to \eqref{eqn:42},
$$
    \dfrac{P(\mathbf W_C\in H_C^1\mid\mathbf  W_A \in L_A^1, \mathbf W_B\in H_B^1)}{P(\mathbf W_C\in H_C^1\mid \mathbf W_A \in L_A^2,\mathbf  W_B\in H_B^1)} = \dfrac{P(\mathbf W_C \in L_C^2\mid \mathbf W_A \in L_A^1, \mathbf W_B\in H_B^1)}{P(\mathbf W_C \in L_C^2\mid \mathbf W_A \in L_A^2, \mathbf W_B\in H_B^1)} = 1.
$$
Due to \eqref{eqn:ind3} and \eqref{eqn:cons3},
it follows directly that 
\begin{multline}
\label{eqn:43}
P(\mathbf W_C\in H_C^2\mid \mathbf W_A \in L_A^1, \mathbf W_B\in H_B^1) \\
= P(\mathbf W_C\in H_C^2\mid \mathbf W_A \in L_A^2, \mathbf W_B\in H_B^1) \\
= P(\mathbf W_C\in H_C^2\mid \mathbf W_A \in H_A^2, \mathbf W_B\in H_B^1) \\
= P(\mathbf W_C\in H_C^2\mid \mathbf W_B\in H_B^1).
\end{multline}
Similarly we have
\begin{multline}
\label{eqn:44}
P(\mathbf W_C\in H_C^1\mid \mathbf W_A \in L_A^1, \mathbf W_B\in H_B^1) \\= P(\mathbf W_C\in H_C^1 \mid \mathbf W_A \in H_A^2, \mathbf W_B\in H_B^1) = P(\mathbf W_C\in H_C^1\mid  \mathbf W_B\in H_B^1).
\end{multline}

We still need to show 
\begin{flalign}
P(\mathbf W_C \in H_C^2 \mid \mathbf W_A \in H_A^1, \mathbf W_B\in H_B^1) &= P(\mathbf W_C \in H_C^2 \mid \mathbf W_A \in L_A^1, \mathbf W_B\in H_B^1), \label{eqn:45}\\
P(\mathbf W_C \in H_C^1 \mid \mathbf W_A \in H_A^1, \mathbf W_B\in H_B^1) &= P(\mathbf W_C \in H_C^1 \mid \mathbf W_A \in L_A^1, \mathbf W_B\in H_B^1).\label{eqn:46}
\end{flalign}

To show \eqref{eqn:45}, note that due to \eqref{eqn:ind3} and \eqref{eqn:cons3},
\begin{small}
\begin{flalign*}
    & \ \ P(\mathbf W_C \in H_C^2\mid \mathbf W_A \in H_A^1, \mathbf W_B\in H_B^1)    \\
    &= \dfrac{P(\mathbf W_C \in H_C^2, \mathbf W_A \in H_A^1, \mathbf W_B\in H_B^1)}{P(\mathbf W_A \in H_A^1, \mathbf W_B\in H_B^1)} \\
    &= \dfrac{ \lambda  P (\mathbf Y_C \in H_C^2, \mathbf Y_C\in L_C^2, \mathbf Y_A \in H_A^1,  \mathbf Y_B\in H_B^1)}{P(\mathbf W_A \in H_A^1, \mathbf W_B\in H_B^1, \mathbf W_C\in L_C^1) + P(\mathbf W_A \in H_A^1, \mathbf W_B\in H_B^1,\mathbf  W_C\in L_C^2)} \\
    &= \dfrac{ \lambda   P(\mathbf Y_C \in H_C^2 \mid  \mathbf Y_B\in H_B^1, \mathbf Y_C\in L_C^2) P(\mathbf Y_A \in H_A^1,  \mathbf Y_B\in H_B^1, \mathbf Y_C\in L_C^2)}{ \dfrac{P(\mathbf W_A \in H_A^1, \mathbf W_B\in H_B^1, \mathbf W_C\in L_C^2){P(\mathbf W_A \in L_A^2, \mathbf W_B\in H_B^1, \mathbf W_C\in L_C^1)}}{P(\mathbf W_A \in L_A^2, \mathbf W_B\in H_B^1, \mathbf W_C\in L_C^2)} + P(\mathbf W_A \in H_A^1, \mathbf W_B\in H_B^1, \mathbf W_C\in L_C^2)} \\
    &= \dfrac{   P(\mathbf Y_C \in H_C^2 \mid  \mathbf Y_B\in H_B^1,\mathbf  Y_C\in L_C^2) P(\mathbf W_A \in H_A^1,  \mathbf W_B\in H_B^1, \mathbf W_C\in L_C^2)}{ \dfrac{P( \mathbf W_A \in H_A^1,  \mathbf W_B\in H_B^1,\mathbf W_C\in L_C^2){P(\mathbf W_A \in L_A^2, \mathbf W_B\in H_B^1, \mathbf W_C\in L_C^1)}}{P(\mathbf W_A \in L_A^2,  \mathbf W_B\in H_B^1,\mathbf  W_C\in L_C^2)} + P(\mathbf W_A \in H_A^1, \mathbf W_B\in H_B^1, \mathbf W_C\in L_C^2)} \\
      &= \dfrac{   P(\mathbf Y_C \in H_C^2 \mid  \mathbf Y_B\in H_B^1, \mathbf Y_C\in L_C^2) }{ \dfrac{{P(\mathbf W_A \in L_A^2, \mathbf W_B\in H_B^1,\mathbf  W_C\in L_C^1)}}{P(\mathbf W_A \in L_A^2,\mathbf  W_B\in H_B^1, \mathbf W_C\in L_C^2)} + 1} 
\end{flalign*}
\end{small}
does not depend on $H_A^1$.

To show \eqref{eqn:46}, note that due to \eqref{eqn:definition3} and \eqref{eqn:cons3},
\begin{small}
\begin{flalign*}
& \ \ P(\mathbf W_C \in H_C^1 \mid \mathbf W_A \in H_A^1, \mathbf W_B\in H_B^1) \\ &= \dfrac{P(\mathbf W_A \in H_A^1, \mathbf W_B\in H_B^1, \mathbf W_C \in H_C^1)}{P(\mathbf W_A \in H_A^1, \mathbf W_B\in H_B^1, \mathbf W_C \in L_C^1) + P(\mathbf W_A \in H_A^1, \mathbf W_B\in H_B^1, \mathbf W_C \in L_C^2)} \\
&= \dfrac{ \dfrac{P(\mathbf W_A \in H_A^1,\mathbf  W_B\in H_B^1, \mathbf W_C \in L_C^2){P(\mathbf W_A \in L_A^2,\mathbf  W_B\in H_B^1, \mathbf W_C \in H_C^1)}}{P(\mathbf W_A \in L_A^2, \mathbf W_B\in H_B^1, \mathbf W_C \in L_C^2)} }{\dfrac{P(\mathbf W_A \in H_A^1, \mathbf W_B\in H_B^1, \mathbf W_C \in L_C^2){P(\mathbf W_A \in L_A^2, \mathbf W_B\in H_B^1, \mathbf W_C \in L_C^1)}}{P(\mathbf W_A \in L_A^2,\mathbf  W_B\in H_B^1, \mathbf W_C \in L_C^2)} + P(\mathbf W_A \in H_A^1, \mathbf W_B\in H_B^1, \mathbf W_C \in L_C^2)} \\
&= \dfrac{ \dfrac{{P(\mathbf W_A \in L_A^2, \mathbf W_B\in H_B^1, \mathbf W_C \in H_C^1)}}{P(\mathbf W_A \in L_A^2, \mathbf W_B\in H_B^1, \mathbf W_C \in L_C^2)} }{\dfrac{{P(\mathbf W_A \in L_A^2, \mathbf W_B\in H_B^1, \mathbf W_C \in L_C^1)}}{P(\mathbf W_A \in L_A^2,\mathbf  W_B\in H_B^1, \mathbf W_C \in L_C^2)} + 1} 
\end{flalign*}
\end{small}
does not depend on $H_A^1$.

From \eqref{eqn:43}, \eqref{eqn:44}, \eqref{eqn:45} and \eqref{eqn:46}, we have
\begin{equation*}
        P( \mathbf W_C \in H_C \mid \mathbf W_A \in H_A, \mathbf W_B \in H_B^1)=P(\mathbf W_C \in H_C \mid  \mathbf W_B \in H_B^1).
\end{equation*}
Besides, because of \eqref{eqn:41}, we can obtain that
\begin{equation*}
        P(\mathbf W_C \in H_C \mid \mathbf W_A \in H_A, \mathbf W_B \in H_B)=P(\mathbf W_C \in H_C \mid  \mathbf W_B \in H_B).
\end{equation*}

 We have hence completed the proof.

\end{proof}

\subsection{Proof of Proposition \ref{prop:eh} in the case where  $B=\emptyset$ }

\begin{proof}
		Use the same notation as before, denote $[0,1]^{|k|}$ as $L_k^1$, $ [0,\infty)^{|k|} \setminus [0,1]^{|k|} $ as $L_k^2$ for $k \in \{A,C\}$. Besides, denote $\{0\}^{|k|}$ as $\bm{0}_k$.
		
		We need to show 
	\begin{equation}
	\label{eqn:claim41}
	\text{the support of} \,\mathbf Y=(\mathbf Y_A,\mathbf Y_C)\,  \text{is}\, \Omega =\left( L_A^2 \times \bm{0}_{C}\right) \cup  \left( \bm{0}_{A}\times L_C^2  \right)\Rightarrow \mathbf Y_A \ind_i \mathbf Y_C,
	\end{equation}
	\begin{equation}
	\label{eqn:claim42}
	\text{the support of} \,\mathbf Y=(\mathbf Y_A,\mathbf Y_C)\,  \text{is}\, \Omega =\left( L_A^2 \times \bm{0}_{C}\right) \cup  \left( \bm{0}_{A}\times L_C^2  \right)\Rightarrow \mathbf Y_A \ind_e \mathbf Y_C,
	\end{equation}
	\begin{equation}
	\label{eqn:claim43}
	\text{the support of} \,\mathbf Y=(\mathbf Y_A,\mathbf Y_C)\,  \text{is}\, \Omega =\left( L_A^2 \times \bm{0}_{C}\right) \cup  \left( \bm{0}_{A}\times L_C^2  \right)\Rightarrow \mathbf Y_A \ind_o \mathbf Y_C.
	\end{equation}
	Under \eqref{eqn:claim41}, \eqref{eqn:claim42} and  \eqref{eqn:claim43}, we can obtain that if the support of $\mathbf Y=(\mathbf Y_A,\mathbf Y_C)$  is $\Omega =\left( L_A^2 \times \bm{0}_{C}\right) \cup  \left( \bm{0}_{A}\times L_C^2  \right)$, then 
	\[
	\mathbf Y_A \ind_o \mathbf Y_C  \Leftrightarrow \mathbf Y_A \ind_i \mathbf Y_C  \Leftrightarrow \mathbf Y_A \ind_e \mathbf Y_C .
	\]
	
	First, we  show \eqref{eqn:claim41}. Since the support of $\mathbf Y=(\mathbf Y_A,\mathbf Y_C)$  is $\Omega =\left( L_A^2 \times \bm{0}_{C}\right) \cup  \left( \bm{0}_{A}\times L_C^2  \right)$, then $S_A\times S_C \subset \left( L_A^2 \times \bm{0}_{C}\right) $ or $S_A\times S_C \subset \left( \bm{0}_{A}\times L_C^2  \right)$. 
	When $S_A\times S_C \subset \left( L_A^2 \times \bm{0}_{C}\right) $, since $\mathbf Y_C=\bm{0}_{C}$ regardless of the value of $\mathbf Y_A$, then $\mathbf Y_A \ind \mathbf Y_C \mid \left(\mathbf Y \in S_A \times S_C\right)$.
	When $S_A\times S_C \subset  \left( \bm{0}_{A}\times L_C^2  \right)$, since $\mathbf Y_A=\bm{0}_{A}$ regardless of the value of $\mathbf Y_C$, then $\mathbf Y_A \ind \mathbf Y_C \mid \left(\mathbf Y \in S_A \times S_C\right)$.
	Thus, based on the definition of $\mathbf Y_A \ind_i \mathbf Y_C$, we have $\mathbf Y_A \ind_i \mathbf Y_C$.

	Second, we show \eqref{eqn:claim42}. Since the support of $\mathbf Y=(\mathbf Y_A,\mathbf Y_C)$  is $\Omega =\left( L_A^2 \times \bm{0}_{C}\right) \cup  \left( \bm{0}_{A}\times L_C^2  \right)$, then when $Y_k>1$ for some $k \in A$, we have $\mathbf Y \in\left( L_A^2 \times \bm{0}_{C}\right)$. Since $\mathbf Y_C=\bm{0}_{C}$ regardless of the value of $\mathbf Y_A$, then $\mathbf Y_A \ind \mathbf Y_C \mid Y_k>1$.  When $Y_k>1$ for some $k\in C$, we have $\mathbf  Y \in \left( \bm{0}_{A}\times L_C^2  \right)$. Since $\mathbf Y_A=\bm{0}_{A}$ regardless of the value of $\mathbf Y_C$, then $\mathbf Y_A \ind \mathbf Y_C \mid Y_k>1$. Thus, based on the definition of $\mathbf Y_A \ind_e \mathbf Y_C$, we have $\mathbf Y_A \ind_e \mathbf Y_C$.
	
	Third, we show \eqref{eqn:claim43}.  Since the support of $\mathbf{Y}=(\mathbf{Y}_A,\mathbf{Y}_C)$  is $\Omega =\left( L_A^2 \times \bm{0}_{C}\right) \cup  \left( \bm{0}_{A}\times L_C^2  \right)$, then 
    \begin{multline*}
	P(\mathbf{Y}_A \in  L_A^2, \mathbf{Y}_C=\bm{0}_C\,\text{or}\, \mathbf{Y}_A=\bm{0}_A, \mathbf{Y}_C \in L_C^2  )\\
	=P(\mathbf{Y}_A \in  L_A^2, \mathbf{Y}_C=\bm{0}_C) + P(\mathbf{Y}_A=\bm{0}_A, \mathbf{Y}_C \in L_C^2)=1.
	\end{multline*}
	Suppose $P(\mathbf{Y}_A \in  L_A^2, \mathbf{Y}_C=\bm{0}_C)=\alpha$ $(0<\alpha<1)$, then $P(\mathbf{Y}_A=\bm{0}_A, \mathbf{Y}_C \in L_C^2)=1-\alpha$. Now we introduce two new random vectors $\mathbf{W}_A$ and $\mathbf W_C$, assume that they have the distributions as follows:
	\[
F_{\mathbf{W}_A}(\mathbf w_A)=
\begin{cases}
	(1-p_A)P(\mathbf U_A\leq \mathbf w_A),& \mathbf w_A \in  L_A^1,\\
	(1-p_A)+p_AP(\mathbf Y_A \leq \mathbf w_A \mid \mathbf Y_A \in L_A^2, \mathbf Y_C=\bm{0}_C), & \mathbf w_A \in L_A^2,
	\end{cases}
	\]
	and 
	\[
F_{\mathbf W_C}(\mathbf w_C)=
\begin{cases}
(1-p_C)P(\mathbf U_C \leq \mathbf w_C), &  \mathbf w_C \in L_C^1,\\
(1-p_C)+p_C P(\mathbf Y_C \leq \mathbf w_C \mid \mathbf Y_A=\bm 0_A, \mathbf Y_C \in L_C^2), & \mathbf w_C \in L_C^2,
\end{cases}
	\]
	where $p_A \in (0,1)$, $p_C \in (0,1)$, $\mathbf U_A$ is a random vector which follows the uniform distribution in $L_A^1$, and $\mathbf U_C$ is a random vector which follows the uniform distribution in $L_C^1$.
	
	To make $\mathbf W_A \ind \mathbf W_C$, let $F(\mathbf w_A,\mathbf w_C)=F_{\mathbf W_A}(\mathbf w_A)F_{\mathbf W_C}(\mathbf w_C)$. Thus, 
	if $\mathbf w_A \in L_A^1, \mathbf w_C \in L_C^1$,
	\[
	 F(\mathbf w_A,\mathbf w_C)= (1-p_A)(1-p_C)P(\mathbf U_A\leq \mathbf w_A)P(\mathbf U_C \leq \mathbf w_C).
	\]
	If $\mathbf w_A \in L_A^1, \mathbf w_C \in L_C^2$,
	\[
	 F(\mathbf w_A, \mathbf w_C)=(1-p_A)P(\mathbf U_A \leq \mathbf w_A)\big\{(1-p_C)+p_C P(\mathbf Y_C \leq \mathbf w_C \mid \mathbf Y_A=\bm 0_A, \mathbf Y_C \in L_C^2)\big\}.
	\]
	If $\mathbf w_A \in L_A^2, \mathbf w_C \in L_C^1$,
	\[
	F(\mathbf w_A,\mathbf w_C)=\big\{(1-p_A)+p_A P(\mathbf Y_A \leq \mathbf w_A \mid \mathbf Y_A \in L_A^2, \mathbf Y_C=\bm 0_C)\big\}(1-p_C)P(\mathbf U_C \leq \mathbf w_C).
	\]
	If $\mathbf w_A \in L_A^2, \mathbf w_C \in L_C^2$, then 
	\begin{multline*}
	F(\mathbf w_A,\mathbf w_C)=\big\{(1-p_A)+p_A P(\mathbf Y_A \leq \mathbf  w_A \mid \mathbf  Y_A \in L_A^2, \mathbf Y_C=\bm 0_C)\big\}\\
	\times \big\{(1-p_C)+p_C P(\mathbf Y_C \leq \mathbf w_C \mid \mathbf Y_A=\bm 0_A, \mathbf Y_C \in L_C^2)\big\}.
	\end{multline*}

If the density function of $F(\mathbf w_A,\mathbf w_C)$ exists, then 
\[
f(\mathbf w_A,\mathbf w_C)=
\begin{cases}
(1-p_A)(1-p_C), & \mathbf w_A \in L_A^1, \mathbf w_C \in L_C^1,\\
\frac{(1-p_A)p_C}{1-\alpha}\frac{\partial P(\mathbf Y_A=\bm 0_A,\mathbf Y_C \leq \mathbf w_C)}{\partial \mathbf w_C},  &   \mathbf w_A \in L_A^1, \mathbf w_C \in L_C^2,\\
\frac{p_A(1-p_C)}{\alpha}\frac{\partial P(\mathbf Y_A\leq \mathbf w_A,\mathbf Y_C = \bm 0_C)}{\partial \mathbf w_A},  &    \mathbf w_A \in L_A^2, \mathbf w_C \in L_C^1,\\
\frac{p_Ap_C}{\alpha(1-\alpha)}\frac{\partial P(\mathbf Y_A\leq \mathbf w_A,\mathbf Y_C = \bm 0_C)}{\partial \mathbf w_A}\frac{\partial P(\mathbf Y_A=\bm 0_A,\mathbf Y_C \leq \mathbf w_C)}{\partial \mathbf w_C}, & \mathbf w_A \in L_A^2, \mathbf w_C \in L_C^2.
\end{cases}
\]

Since the support of $\mathbf Y = (\mathbf Y_A,\mathbf Y_C) $ is $\Omega =\left( L_A^2 \times \bm{0}_{C}\right) \cup  \left( \bm{0}_{A}\times L_C^2  \right)$, then $\partial P(\mathbf Y_A=\bm 0_A,\mathbf Y_C \leq \mathbf w_C)/\partial \mathbf w_C \neq 0$ for $\mathbf w_C  \in L_C^2$ and $\partial P(\mathbf Y_A\leq \mathbf w_A,\mathbf Y_C = \bm 0_C)/\partial \mathbf w_A\neq 0$ for $\mathbf w_A \in L_A^2$. Thus, from $f(\mathbf w_1,\mathbf w_2)$, we have the support of $\mathbf W=(\mathbf W_A,\mathbf W_C)$ is $[0,\infty)^{|A|+|C|}$, which contains  $\Omega$.  Based on the definition of conditional density in \citet[\S 7.4]{shorack2017probability}, we have that the conditional density at $(\mathbf w_A,\bm 0_C)$, $(\mathbf w_A \in L_A^2)$ and $(\bm 0_A,\mathbf w_C)$, $(\mathbf w_C \in L_C^2)$ given that all the points are on the curve of $\Omega$ are
\[
\begin{aligned}
&\frac{f(\mathbf w_A,\bm 0_C)}{ \int_{L_A^2} f(\mathbf w_A,\bm 0_C) {\rm d} \mathbf w_A  + \int_{L_C^2} f(\bm 0_A, \mathbf w_C) {\rm d} \mathbf w_C  }\\
=&\frac{ \frac{p_A(1-p_C)}{\alpha}\frac{\partial P(\mathbf Y_A\leq \mathbf w_A,\mathbf Y_C = \bm 0_C)}{\partial \mathbf w_A}}{ \frac{p_A(1-p_C)}{\alpha}P(\mathbf Y_A \in L_A^2,\mathbf Y_C=\bm 0_C)+\frac{(1-p_A)p_C}{1-\alpha}P(\mathbf Y_A=\bm 0_A,\mathbf Y_C\in L_C^2)}\\
=&\frac{p_A(1-p_C)\frac{\partial P(\mathbf Y_A\leq \mathbf w_A, \mathbf Y_C=\bm 0_C )}{\partial \mathbf w_A}
}{\alpha\{p_A(1-p_C)+ (1-p_A)p_C \}   },
\end{aligned}
\]
and
\[
\begin{aligned}
&\frac{f(\bm 0_A,\mathbf w_C)}{ \int_{L_A^2} f(\mathbf w_A,\bm 0_C) {\rm d} \mathbf w_A  + \int_{L_C^2} f(\bm 0_A, \mathbf w_C) {\rm d} \mathbf w_C  }\\
=&\frac{\frac{(1-p_A)p_C}{1-\alpha}\frac{\partial P(\mathbf Y_A=\bm 0_A,\mathbf Y_C \leq \mathbf w_C)}{\partial \mathbf w_C}}{ \frac{p_A(1-p_C)}{\alpha}P(\mathbf Y_A\in L_A^2,\mathbf Y_C=\bm 0_C)+\frac{(1-p_A)p_C}{1-\alpha}P(\mathbf Y_A=\bm 0_A,\mathbf Y_C \in L_C^2)}\\
=&\frac{(1-p_A)p_C\frac{\partial P(\mathbf Y_A=\bm 0_A,\mathbf Y_C\leq \mathbf w_C)}{\partial \mathbf w_C} }{(1-\alpha)\{p_A(1-p_C)+ (1-p_A)p_C\}   }.
\end{aligned}
\]

Then we can get the value of $P(\mathbf W_A\leq \mathbf w_A, \mathbf W_C\leq \mathbf w_C\mid \mathbf W_A\in L_A^2, \mathbf  W_C=\bm 0_C \, \text{or}\, \mathbf W_A=\bm 0_A, \mathbf W_C\in L_C^2)$. When  $\mathbf w_A \in L_A^1, \mathbf w_C \in L_C^1$, 
\[
P(\mathbf W_A\leq \mathbf w_A, \mathbf W_C\leq \mathbf w_C\mid \mathbf W_A\in L_A^2, \mathbf  W_C=\bm 0_C \, \text{or}\, \mathbf W_A=\bm 0_A, \mathbf W_C\in L_C^2)=0=P(\mathbf Y_A \leq \mathbf w_A, \mathbf Y_C \leq \mathbf w_C).
\]
When $\mathbf w_A \in L_A^1, \mathbf w_C \in L_C^2$, 
\[
\begin{aligned}
&P(\mathbf W_A\leq \mathbf w_A, \mathbf W_C\leq \mathbf w_C\mid \mathbf W_A\in L_A^2, \mathbf  W_C=\bm 0_C \, \text{or}\, \mathbf W_A=\bm 0_A, \mathbf W_C\in L_C^2)\\
=&\frac{ \int _{[\bm{0}_C,\mathbf w_C]\setminus[0,1]^{|C|} } f(\bm 0_A,\mathbf t) {\rm d} \mathbf t  }{\int _{L_A^2} f(\mathbf w_A,\bm 0_C){\rm d} \mathbf w_A+\int_{L_C^2} f(\bm 0_A,\mathbf w_C){\rm d} \mathbf w_C}\\
=&\frac{(1-p_A)p_CP(\mathbf Y_A=\bm 0_A,\mathbf Y_C \leq \mathbf w_C)  }{(1-\alpha)\{p_A(1-p_C)+(1-p_A)p_C  \}}\\
=&\frac{(1-p_A)p_CP(\mathbf Y_A\leq \mathbf w_A,\mathbf Y_C \leq \mathbf w_C)  }{(1-\alpha)\{p_A(1-p_C)+(1-p_A)p_C  \}},
\end{aligned}
\]
where the last equation is due to the support of $\mathbf{Y}=(\mathbf{Y}_A,\mathbf{Y}_C)$  being $\Omega =\left( L_A^2 \times \bm{0}_{C}\right) \cup  \left( \bm{0}_{A}\times L_C^2  \right)$.
When $\mathbf w_A \in L_A^2, \mathbf w_C \in L_C^1$, 
\[
\begin{aligned}
&P(\mathbf W_A\leq \mathbf w_A, \mathbf W_C\leq \mathbf w_C\mid \mathbf W_A\in L_A^2, \mathbf  W_C=\bm 0_C \, \text{or}\, \mathbf W_A=\bm 0_A, \mathbf W_C\in L_C^2)\\
=&\frac{ \int _{[\bm{0}_A,\mathbf w_A]\setminus[0,1]^{|A|} }  f(\mathbf t,\bm 0_C) {\rm d} \mathbf t      }{\int _{L_A^2} f(\mathbf w_A,\bm 0_C){\rm d} \mathbf w_A+\int_{L_C^2} f(\bm 0_A,\mathbf w_C){\rm d} \mathbf w_C}\\
=&\frac{p_A(1-p_C)P(\mathbf Y_A \leq \mathbf w_A,\mathbf Y_C = \bm 0_C)  }{\alpha\{p_A(1-p_C)+(1-p_A)p_C \}}\\
=&\frac{p_A(1-p_C)P(\mathbf Y_A\leq \mathbf w_A,\mathbf Y_C \leq \mathbf w_C)  }{\alpha\{p_A(1-p_C)+(1-p_A)p_C\}}.
\end{aligned}
\]
When $\mathbf w_A \in L_A^2, \mathbf w_C \in L_C^2$, due to $P(\mathbf Y_A \leq \mathbf w_A, \mathbf Y_C \leq \mathbf w_C)=P(\mathbf Y_A \leq \mathbf w_A, \mathbf Y_C = \bm 0_C)+P(\mathbf Y_A = \bm 0_A, \mathbf Y_C \leq \mathbf w_C)$,
\[
\begin{aligned}
&P(\mathbf W_A\leq \mathbf w_A, \mathbf W_C\leq \mathbf w_C\mid \mathbf W_A\in L_A^2, \mathbf  W_C=\bm 0_C \, \text{or}\, \mathbf W_A=\bm 0_A, \mathbf W_C\in L_C^2)\\
=&\frac{\int _{[\bm{0}_A,\mathbf w_A]\setminus[0,1]^{|A|} }  f(\mathbf t,\bm 0_C) {\rm d} \mathbf t }{\int _{L_A^2} f(\mathbf w_A,\bm 0_C){\rm d} \mathbf w_A+\int_{L_C^2} f(\bm 0_A,\mathbf w_C){\rm d} \mathbf w_C}+\frac{\int _{[\bm{0}_C,\mathbf w_C]\setminus[0,1]^{|C|} } f(\bm 0_A,\mathbf t) {\rm d} \mathbf t  }{\int _{L_A^2} f(\mathbf w_A,\bm 0_C){\rm d} \mathbf w_A+\int_{L_C^2} f(\bm 0_A,\mathbf w_C){\rm d} \mathbf w_C}\\
=&\frac{p_A(1-p_C)P(\mathbf Y_A \leq \mathbf w_A,\mathbf Y_C = \bm 0_C)  }{\alpha\{p_A(1-p_C)+(1-p_A)p_C \}}+\frac{(1-p_A)p_CP(\mathbf Y_A=\bm 0_A,\mathbf Y_C \leq \mathbf w_C)  }{(1-\alpha)\{p_A(1-p_C)+(1-p_A)p_C  \}}\\
=&\frac{\{(1-\alpha)p_A(1-p_C)-\alpha (1-p_A)p_C\}P(\mathbf Y_A\leq \mathbf w_A,\mathbf Y_C = \bm 0_C)  }{\alpha(1-\alpha)\{p_A(1-p_C)+(1-p_A)p_C  \}}+\frac{(1-p_A)p_CP(\mathbf Y_A\leq \mathbf w_A,\mathbf Y_C \leq \mathbf w_C)  }{(1-\alpha)\{p_A(1-p_C)+(1-p_A)p_C  \}}.
\end{aligned}
\]

In summary,
\begin{equation}
\label{eqn:50}
\begin{aligned}
&P(\mathbf W_A\leq \mathbf w_A, \mathbf W_C\leq \mathbf w_C\mid \mathbf W_A\in L_A^2, \mathbf  W_C=\bm 0_C \, \text{or}\, \mathbf W_A=\bm 0_A, \mathbf W_C\in L_C^2)=\\
&\begin{cases}
0, &\mathbf w_A \in L_A^1, \mathbf w_C \in L_C^1,\\
\frac{(1-p_A)p_CP(\mathbf Y_A\leq \mathbf w_A,\mathbf Y_C \leq \mathbf w_C)  }{(1-\alpha)\{p_A(1-p_C)+(1-p_A)p_C\}}, & \mathbf w_A \in L_A^1, \mathbf w_C \in L_C^2,\\
\frac{p_A(1-p_C)P(\mathbf Y_A\leq \mathbf w_A,\mathbf Y_C \leq \mathbf w_C)  }{\alpha\{p_A(1-p_C)+(1-p_A)p_C\}}, & \mathbf w_A \in L_A^2, \mathbf w_C \in L_C^1,\\
 \frac{\{(1-\alpha)p_A(1-p_C)-\alpha (1-p_A)p_C\}P(\mathbf Y_A\leq \mathbf w_A,\mathbf Y_C = \bm 0_C)  }{\alpha(1-\alpha)\{p_A(1-p_C)+(1-p_A)p_C  \}}+\frac{(1-p_A)p_CP(\mathbf Y_A\leq \mathbf w_A,\mathbf Y_C \leq \mathbf w_C)  }{(1-\alpha)\{p_A(1-p_C)+(1-p_A)p_C  \}},& \mathbf w_A \in L_A^2, \mathbf w_C \in L_C^2.
\end{cases} 
\end{aligned}
\end{equation}
From \eqref{eqn:50},  $P(\mathbf W_A\leq \mathbf w_A, \mathbf W_C\leq \mathbf w_C\mid \mathbf W_A\in L_A^2, \mathbf  W_C=\bm 0_C \, \text{or}\, \mathbf W_A=\bm 0_A, \mathbf W_C\in L_C^2)=P(\mathbf Y_A\leq \mathbf w_A,\mathbf Y_C \leq \mathbf w_C)$,  if and only if 
\[
\begin{cases}
\frac{(1-p_A)p_C}{(1-\alpha)\{p_A(1-p_C)+(1-p_A)p_C  \}  }=1,\\
\frac{p_A(1-p_C)}{\alpha\{p_A(1-p_C)+(1-p_A)p_C  \}  }=1,\\
\frac{\{(1-\alpha)p_A(1-p_C) - \alpha(1-p_A)p_C \} }{\alpha(1-\alpha)\{ p_A(1-p_C)+(1-p_A)p_C  \} }=0,
\end{cases}
\]
or equivalently,
\begin{equation}
\label{eqn:51}
\alpha=\frac{p_A(1-p_C)}{p_A(1-p_C)+(1-p_A)p_C}.
\end{equation}
We have hence shown that under condition \eqref{eqn:51}, $(\mathbf W \mid \mathbf W\in \Omega) \,{\buildrel d \over =}\, \mathbf Y$ if the density function of $F(\mathbf w_A,\mathbf w_C)$ exists.

If the density function of  $F(\mathbf w_A,\mathbf w_C)$  does not exist, from the definition of conditional probability in \citet[\S 7.4]{shorack2017probability}, we know that \eqref{eqn:50} still satisfies the condition for conditional probability. Thus, we can define $P(\mathbf W_A\leq \mathbf w_A, \mathbf W_C\leq \mathbf w_C\mid \mathbf W_A\in L_A^2, \mathbf  W_C=\bm 0_C \, \text{or}\, \mathbf W_A=\bm 0_A, \mathbf W_C\in L_C^2)$ as that in \eqref{eqn:50}. Under condition \eqref{eqn:51}, we also have $P(\mathbf W_A\leq \mathbf w_A, \mathbf W_C\leq \mathbf w_C\mid \mathbf W_A\in L_A^2, \mathbf  W_C=\bm 0_C \, \text{or}\, \mathbf W_A=\bm 0_A, \mathbf W_C\in L_C^2)=P(\mathbf Y_A\leq \mathbf w_A,\mathbf Y_C \leq \mathbf w_C)$. Then $(\mathbf W \mid \mathbf W\in \Omega) \,{\buildrel d \over =}\, \mathbf Y$.

In summary, based on the definition of $\mathbf Y_A \ind_o \mathbf Y_C$, we have $\mathbf Y_A \ind_o \mathbf Y_C$.

\end{proof}

\subsection{Proof of Proposition \ref{prop:eh} for a general $B$}

\begin{proof}
	We need to show 
	\begin{multline}
	\label{eqn:claim61}
	\text{the support of} \,\mathbf Y=(\mathbf Y_A,\mathbf Y_B,\mathbf Y_C)\,  \text{is}\, \Omega =\left( L_A^2 \times \bm{0}_{B}\times \bm{0}_{C} \right) \cup  \left( \bm{0}_{A}\times L_B^2\times  \bm{0}_{C} \right) \cup  \left( \bm{0}_{A}\times \bm{0}_{B}\times L_C^2  \right)\\
	\Rightarrow \mathbf Y_A \ind_i \mathbf Y_C \mid  \mathbf Y_B,
	\end{multline}
	\begin{multline}
	\label{eqn:claim62}
	\text{the support of} \,\mathbf Y=(\mathbf Y_A,\mathbf Y_B,\mathbf Y_C)\,  \text{is}\, \Omega =\left( L_A^2 \times \bm{0}_{B}\times \bm{0}_{C} \right) \cup  \left( \bm{0}_{A}\times L_B^2\times  \bm{0}_{C} \right) \cup  \left( \bm{0}_{A}\times \bm{0}_{B}\times L_C^2  \right)\\
\Rightarrow \mathbf Y_A \ind_e \mathbf Y_C \mid  \mathbf Y_B,
	\end{multline}
	\begin{multline}
	\label{eqn:claim63}
	\text{the support of} \,\mathbf Y=(\mathbf Y_A,\mathbf Y_B,\mathbf Y_C)\,  \text{is}\, \Omega =\left( L_A^2 \times \bm{0}_{B}\times \bm{0}_{C} \right) \cup  \left( \bm{0}_{A}\times L_B^2\times  \bm{0}_{C} \right) \cup  \left( \bm{0}_{A}\times \bm{0}_{B}\times L_C^2  \right)\\
\Rightarrow \mathbf Y_A \ind_o \mathbf Y_C \mid  \mathbf Y_B.
	\end{multline}
	Under \eqref{eqn:claim61}, \eqref{eqn:claim62} and  \eqref{eqn:claim63}, we can have if the support of $\mathbf Y=(\mathbf Y_A,\mathbf Y_B,\mathbf Y_C)$  is $\Omega =\left( L_A^2 \times \bm{0}_{B}\times \bm{0}_{C} \right) \cup  \left( \bm{0}_{A}\times L_B^2\times  \bm{0}_{C} \right) \cup  \left( \bm{0}_{A}\times \bm{0}_{B}\times L_C^2  \right)$, then 
	\[
	\mathbf Y_A \ind_o \mathbf Y_C \mid  \mathbf Y_B \Leftrightarrow \mathbf Y_A \ind_i \mathbf Y_C \mid  \mathbf Y_B\Leftrightarrow \mathbf Y_A \ind_e \mathbf Y_C \mid  \mathbf Y_B.
	\]
	
	First, we  show \eqref{eqn:claim61}. Since the support of $\mathbf Y=(\mathbf Y_A,\mathbf Y_B,\mathbf Y_C)$  is $\Omega =\left( L_A^2 \times \bm{0}_{B}\times \bm{0}_{C} \right) \cup  \left( \bm{0}_{A}\times L_B^2\times  \bm{0}_{C} \right) \\
	\cup  \left( \bm{0}_{A}\times \bm{0}_{B}\times L_C^2  \right)$, then $S_A\times S_B\times S_C \subset \left( L_A^2 \times \bm{0}_{B}\times \bm{0}_{C} \right)$, $\left( \bm{0}_{A}\times L_B^2\times  \bm{0}_{C} \right)$ or $\left( \bm{0}_{A}\times \bm{0}_{B}\times L_C^2  \right)$. 
	When $S_A\times S_B \times S_C \subset\left( L_A^2 \times \bm{0}_{B}\times \bm{0}_{C} \right) $, $\mathbf Y_B=\bm{0}_{B}$. Given $\mathbf Y_B=\bm{0}_{B}$, since $\mathbf Y_C=\bm{0}_{C}$ regardless of the value of $\mathbf Y_A$, then $\mathbf Y_A \ind \mathbf Y_C \mid (\mathbf Y_B, \mathbf Y \in S_A \times S_B \times S_C)$.
	When $S_A\times S_B\times S_C \subset\left( \bm{0}_{A}\times L_B^2\times  \bm{0}_{C} \right)$, $\mathbf Y_B \in L_B^2$. Given $\mathbf Y_B=\mathbf y_B$, since $\mathbf Y_A=\bm{0}_{A}$ and  $\mathbf Y_C=\bm{0}_{C}$,  then $\mathbf Y_A \ind \mathbf Y_C \mid (\mathbf Y_B, \mathbf Y \in S_A \times S_B \times S_C)$. When $S_A\times S_B\times S_C \subset \left( \bm{0}_{A}\times \bm{0}_{B}\times L_C^2  \right)$, $\mathbf Y_B=\bm{0}_{B}$. Given $\mathbf Y_B=\bm{0}_{B}$, since $\mathbf Y_A=\bm{0}_{A}$ regardless of the value of $\mathbf Y_C$, then $\mathbf Y_A \ind \mathbf Y_C \mid (\mathbf Y_B, \mathbf Y \in S_A \times S_B \times S_C)$.
	Thus, based on the definition of $\mathbf Y_A \ind_i \mathbf Y_C \mid  \mathbf Y_B$, we have $\mathbf Y_A \ind_i \mathbf Y_C \mid  \mathbf Y_B$.

	Second, we show \eqref{eqn:claim62}. Since the support of $\mathbf Y=(\mathbf Y_A,\mathbf Y_B,\mathbf Y_C)$  is $\Omega =\left( L_A^2 \times \bm{0}_{B}\times \bm{0}_{C} \right) \cup  \left( \bm{0}_{A}\times L_B^2\times  \bm{0}_{C} \right)\\ \cup  \left( \bm{0}_{A}\times \bm{0}_{B}\times L_C^2  \right)$, then when $Y_k>1$ for some $k\in A$, we have $\mathbf Y \in \left( L_A^2 \times \bm{0}_{B}\times \bm{0}_{C} \right)$ and $\mathbf{Y}_B=\bm{0}_{B}$. Given $\mathbf Y_B=\bm{0}_{B}$ and $Y_k>1$, since $\mathbf{Y}_C=\bm{0}_{C}$ regardless of the value of $\mathbf Y_A$, then $\mathbf Y_A \ind \mathbf Y_C \mid  (\mathbf Y_B, Y_k>1)$.  When $Y_k>1$ for some $k\in B$, we have $\mathbf Y \in   \left( \bm{0}_{A}\times L_B^2\times  \bm{0}_{C} \right)$ and $\mathbf Y_B \in  L_B^2$. Given $\mathbf Y_B=\mathbf y_B$ and $Y_k>1$, since $\mathbf Y_A=\bm{0}_{A}$ and  $\mathbf Y_C=\bm{0}_{C}$, then $\mathbf Y_A \ind \mathbf Y_C \mid  (\mathbf Y_B, Y_k>1)$. When $Y_k>1$ for some $k\in C$, we have $\mathbf Y \in    \left( \bm{0}_{A}\times \bm{0}_{B}\times L_C^2  \right)$ and $\mathbf Y_B =\bm{0}_{B}$. Given $\mathbf Y_B=\bm{0}_{B}$ and $Y_k>1$, since $\mathbf Y_A=\bm{0}_{A}$ regardless of the value of $\mathbf Y_C$, then $\mathbf Y_A \ind \mathbf Y_C \mid  (\mathbf Y_B, Y_k>1)$.
	Thus, based on the definition of $\mathbf Y_A \ind_e \mathbf Y_C \mid  \mathbf Y_B$, we have $\mathbf Y_A \ind_e \mathbf Y_C \mid  \mathbf Y_B$.
	
	Third, we show \eqref{eqn:claim63}. 
	Since the support of $\mathbf Y=(\mathbf Y_A,\mathbf Y_B,\mathbf Y_C)$  is $\Omega =\left( L_A^2 \times \bm{0}_{B}\times \bm{0}_{C} \right) \cup  \left( \bm{0}_{A}\times L_B^2\times  \bm{0}_{C} \right)\\ \cup  \left( \bm{0}_{A}\times \bm{0}_{B}\times L_C^2  \right)$, then 
\begin{multline*}
P(\mathbf Y_A\in L_A^2,\mathbf Y_B=\bm{0}_{B},\mathbf Y_C=\bm{0}_{C} \,\text{or}\, \mathbf Y_A=\bm{0}_{A}, \mathbf Y_B\in L_B^2, \mathbf Y_C=\bm{0}_{C}\,\text{or}\,\mathbf Y_A=\bm{0}_{A}, \mathbf Y_B=\bm{0}_{B}, \mathbf Y_C\in L_C^2)\\=P(\mathbf Y_A\in L_A^2,\mathbf Y_B=\bm{0}_{B},\mathbf Y_C=\bm{0}_{C}) + P(\mathbf Y_A=\bm{0}_{A}, \mathbf Y_B\in L_B^2, \mathbf Y_C=\bm{0}_{C})\\
+ P(\mathbf Y_A=\bm{0}_{A}, \mathbf Y_B=\bm{0}_{B}, \mathbf Y_C\in L_C^2)=1.
\end{multline*}
Suppose that $P(\mathbf Y_A\in L_A^2,\mathbf Y_B=\bm{0}_{B},\mathbf Y_C=\bm{0}_{C})=\alpha_1$ and $P(\mathbf Y_A=\bm{0}_{A}, \mathbf Y_B\in L_B^2, \mathbf Y_C=\bm{0}_{C})=\alpha_2$, then $P(\mathbf Y_A=\bm{0}_{A}, \mathbf Y_B=\bm{0}_{B}, \mathbf Y_C\in L_C^2)=1-\alpha_1-\alpha_2$, where $\alpha_1$ and $\alpha_2$ satisfy the conditions of $\alpha_1>0, \alpha_2>0$ and $0<\alpha_1+\alpha_2<1$. 
	Now we introduce three new random vectors $\mathbf W_A$, $\mathbf W_B$ and $\mathbf W_C$, assume that they have the distributions as follows:
	\[
F_{\mathbf{W}_A}(\mathbf w_A)=
\begin{cases}
(1-p_A)P(\mathbf U_A\leq \mathbf w_A),& \mathbf w_A \in  L_A^1,\\
(1-p_A)+p_AP(\mathbf Y_A \leq \mathbf w_A \mid \mathbf Y_A \in L_A^2, \mathbf Y_B=\bm{0}_B, \mathbf Y_C=\bm{0}_C), & \mathbf w_A \in L_A^2,
\end{cases}
\]
\[
F_{\mathbf W_B}(\mathbf w_B)=
\begin{cases}
(1-p_B)P(\mathbf U_B \leq \mathbf w_B), &  \mathbf w_B \in L_B^1,\\
(1-p_B)+p_B P(\mathbf Y_B \leq \mathbf w_B \mid \mathbf Y_A=\bm{0}_A, \mathbf Y_B \in L_B^2, \mathbf Y_C=\bm{0}_C), & \mathbf w_B \in L_B^2,
\end{cases}
\]
and 
\[
F_{\mathbf W_C}(\mathbf w_C)=
\begin{cases}
(1-p_C)P(\mathbf U_C \leq \mathbf w_C), &  \mathbf w_C \in L_C^1,\\
(1-p_C)+p_C P(\mathbf Y_C \leq \mathbf w_C \mid \mathbf Y_A=\bm{0}_A, \mathbf Y_B=\bm{0}_B,\mathbf Y_C \in L_C^2), & \mathbf w_C \in L_C^2,
\end{cases}
\]
where $p_A \in (0,1)$, $p_B \in (0,1)$ $p_C \in (0,1)$, and  $\mathbf U_A$, $\mathbf U_B$, $\mathbf U_C$ are random vectors which follow the uniform distributions in $L_A^1$,  $L_B^1$ and  $L_C^1$, respectively.

To make $\mathbf W_A,\mathbf W_B,\mathbf W_C$ independent, let $F(\mathbf w_A,\mathbf w_B,\mathbf w_C)=F_{\mathbf W_A}(\mathbf w_A)F_{\mathbf W_B}(\mathbf w_B)F_{\mathbf W_C}(\mathbf w_C)$.  Furthermore, we have $\mathbf W_A \ind \mathbf W_C \mid \mathbf W_B$. Besides, 
	if $\mathbf w_A \in L_A^1, \mathbf w_B \in L_B^1, \mathbf w_C \in L_C^1$,
	\[
	 F(\mathbf w_A,\mathbf w_B,\mathbf w_C)= (1-p_A)(1-p_B)(1-p_C)P(\mathbf U_A\leq \mathbf w_A)P(\mathbf U_B\leq \mathbf w_B)P(\mathbf U_C \leq \mathbf w_C).
	\]
	If $\mathbf w_A \in L_A^1, \mathbf w_B \in L_B^1, \mathbf w_C \in L_C^2$,
	\begin{multline*}
	 F(\mathbf w_A, \mathbf w_B, \mathbf w_C)=(1-p_A)(1-p_B)P(\mathbf U_A \leq \mathbf w_A)P(\mathbf U_B \leq \mathbf w_B)\\
	 \times \big\{(1-p_C)+p_C P(\mathbf Y_C \leq \mathbf w_C \mid \mathbf Y_A=\bm 0_A, \mathbf Y_B=\bm 0_B, \mathbf Y_C \in L_C^2)\big\}.	    
	\end{multline*}
	If $\mathbf w_A \in L_A^1, \mathbf w_B \in L_B^2, \mathbf w_C \in L_C^1$,
	\begin{multline*}
	 F(\mathbf w_A, \mathbf w_B, \mathbf w_C)=(1-p_A)(1-p_C)P(\mathbf U_A \leq \mathbf w_A)P(\mathbf U_C \leq \mathbf w_C)\\
	 \times \big\{(1-p_B)+p_B P(\mathbf Y_B \leq \mathbf w_B \mid \mathbf Y_A=\bm 0_A,  \mathbf Y_B \in L_B^2, \mathbf Y_C=\bm 0_C)\big\}.	    
	\end{multline*}
	If $\mathbf w_A \in L_A^1, \mathbf w_B \in L_B^2, \mathbf w_C \in L_C^2$,
	\begin{multline*}
	 F(\mathbf w_A, \mathbf w_B, \mathbf w_C)=(1-p_A)P(\mathbf U_A \leq \mathbf w_A) \big\{(1-p_B)+p_B P(\mathbf Y_B \leq \mathbf w_B \mid \mathbf Y_A=\bm 0_A,  \mathbf Y_B \in L_B^2, \mathbf Y_C=\bm 0_C)\big\}\\
	 \times \big\{(1-p_C)+p_C P(\mathbf Y_C \leq \mathbf w_C \mid \mathbf Y_A=\bm 0_A, \mathbf Y_B=\bm 0_B, \mathbf Y_C \in L_C^2)\big\}.	    
	\end{multline*}
	If $\mathbf w_A \in L_A^2, \mathbf w_B \in L_B^1, \mathbf w_C \in L_C^1$,
	\begin{multline*}
	 F(\mathbf w_A, \mathbf w_B, \mathbf w_C)=(1-p_B)(1-p_C)P(\mathbf U_B \leq \mathbf w_B)P(\mathbf U_C \leq \mathbf w_C)\\
	 \times \big\{(1-p_A)+p_AP(\mathbf Y_A \leq \mathbf w_A \mid \mathbf Y_A \in L_A^2, \mathbf Y_B=\bm{0}_B, \mathbf Y_C=\bm{0}_C)\big\}.	    
	\end{multline*}
		If $\mathbf w_A \in L_A^2, \mathbf w_B \in L_B^1, \mathbf w_C \in L_C^2$,
	\begin{multline*}
	 F(\mathbf w_A, \mathbf w_B, \mathbf w_C)=(1-p_B)P(\mathbf U_B \leq \mathbf w_B) \big\{(1-p_A)+p_AP(\mathbf Y_A \leq \mathbf w_A \mid \mathbf Y_A \in L_A^2, \mathbf Y_B=\bm{0}_B, \mathbf Y_C=\bm{0}_C)\big\}\\
	 \times \big\{(1-p_C)+p_C P(\mathbf Y_C \leq \mathbf w_C \mid \mathbf Y_A=\bm 0_A, \mathbf Y_B=\bm 0_B, \mathbf Y_C \in L_C^2)\big\}.	    
	\end{multline*}
	If $\mathbf w_A \in L_A^2, \mathbf w_B \in L_B^2, \mathbf w_C \in L_C^1$,
	\begin{multline*}
	 F(\mathbf w_A, \mathbf w_B, \mathbf w_C)=(1-p_C)P(\mathbf U_C \leq \mathbf w_C) \big\{(1-p_A)+p_AP(\mathbf Y_A \leq \mathbf w_A \mid \mathbf Y_A \in L_A^2, \mathbf Y_B=\bm{0}_B, \mathbf Y_C=\bm{0}_C)\big\}\\
	 \times  \big\{(1-p_B)+p_B P(\mathbf Y_B \leq \mathbf w_B \mid \mathbf Y_A=\bm 0_A,  \mathbf Y_B \in L_B^2, \mathbf Y_C=\bm 0_C)\big\}.	    
	\end{multline*}
		If $\mathbf w_A \in L_A^2, \mathbf w_B \in L_B^2, \mathbf w_C \in L_C^2$,
	\begin{multline*}
	 F(\mathbf w_A, \mathbf w_B, \mathbf w_C)= \big\{(1-p_A)+p_AP(\mathbf Y_A \leq \mathbf w_A \mid \mathbf Y_A \in L_A^2, \mathbf Y_B=\bm{0}_B, \mathbf Y_C=\bm{0}_C)\big\}\\
	 \times  \big\{(1-p_B)+p_B P(\mathbf Y_B \leq \mathbf w_B \mid \mathbf Y_A=\bm 0_A,  \mathbf Y_B \in L_B^2, \mathbf Y_C=\bm 0_C)\big\}\\
	 \times \big\{(1-p_C)+p_C P(\mathbf Y_C \leq \mathbf w_C \mid \mathbf Y_A=\bm 0_A, \mathbf Y_B=\bm 0_B, \mathbf Y_C \in L_C^2)\big\}.	    
	\end{multline*}

	If the density function of $F(\mathbf w_A,\mathbf w_B,\mathbf w_C)$ exists, then by calculation
	\[
\begin{aligned}
&f(\mathbf w_A,\mathbf w_B,\mathbf w_C)=\\
&\begin{cases}
&(1-p_A)(1-p_B)(1-p_C), \hspace{17.3em} \mathbf w_A \in L_A^1,\mathbf w_B \in L_B^1, \mathbf w_C \in L_C^1,\\
&\frac{(1-p_A)(1-p_B)p_C}{1-\alpha_1-\alpha_2}\frac{\partial P(\mathbf Y_A=\bm 0_A,\mathbf Y_B=\bm 0_B,\mathbf Y_C \leq \mathbf w_C)}{\partial \mathbf w_C},  \hspace{10.6em}    \mathbf w_A \in L_A^1,\mathbf w_B \in L_B^1, \mathbf w_C \in L_C^2,\\
&\frac{(1-p_A)p_B(1-p_C)}{\alpha_2}\frac{\partial P(\mathbf Y_A=\bm 0_A,\mathbf Y_B \leq \mathbf w_B,\mathbf Y_C=\bm 0_C)}{\partial \mathbf w_B},  \hspace{10.7em}    \mathbf w_A \in L_A^1,\mathbf w_B \in L_B^2, \mathbf w_C \in L_C^1,\\
&\frac{(1-p_A)p_B p_C}{\alpha_2(1-\alpha_1-\alpha_2)}\frac{\partial P(\mathbf Y_A=\bm 0_A,\mathbf Y_B\leq \mathbf w_B,\mathbf Y_C=\bm 0_C)}{\partial \mathbf w_B}\frac{\partial P(\mathbf Y_A=\bm 0_A,\mathbf Y_B=\bm 0_B,\mathbf Y_C\leq \mathbf w_C)}{\partial \mathbf w_C},  \hspace{0.7em} \mathbf w_A \in L_A^1,\mathbf w_B \in L_B^2, \mathbf w_C \in L_C^2,\\
&\frac{p_A(1-p_B)(1-p_C)}{\alpha_1}\frac{\partial P(\mathbf Y_A\leq \mathbf w_A,\mathbf Y_B=\bm 0_B,\mathbf Y_C=\bm 0_C)}{\partial \mathbf w_A}, \hspace{10.7em} \mathbf w_A \in L_A^2,\mathbf w_B \in L_B^1, \mathbf w_C \in L_C^1,\\
&\frac{p_A(1-p_B)p_C}{\alpha_1(1-\alpha_1-\alpha_2)}\frac{\partial P(\mathbf Y_A\leq \mathbf w_A,\mathbf Y_B=\bm 0_B,\mathbf Y_C=\bm 0_C)}{\partial \mathbf w_A}\frac{\partial P(\mathbf Y_A=\bm 0_A,\mathbf Y_B=\bm 0_B,\mathbf Y_C\leq \mathbf w_C)}{\partial \mathbf w_C}, \hspace{0.7em} \mathbf w_A \in L_A^2,\mathbf w_B \in L_B^1, \mathbf w_C \in L_C^2,\\
&\frac{p_A p_B(1-p_C)}{\alpha_1\alpha_2}\frac{\partial P(\mathbf Y_A\leq \mathbf w_A,\mathbf Y_B=\bm 0_B,\mathbf Y_C=\bm 0_C)}{\partial \mathbf w_A}\frac{\partial P(\mathbf Y_A=\bm 0_A,\mathbf Y_B\leq \mathbf w_B,\mathbf Y_C=\bm 0_C)}{\partial \mathbf w_B}, \hspace{1em} \mathbf w_A \in L_A^2,\mathbf w_B \in L_B^2, \mathbf w_C \in L_C^1,\\
&\frac{p_A p_B p_C}{\alpha_1\alpha_2(1-\alpha_1-\alpha_2)}\frac{\partial P(\mathbf Y_A\leq \mathbf w_A,\mathbf Y_B=\bm 0_B,\mathbf Y_C=\bm 0_C)}{\partial \mathbf w_A}\\
&\hspace{2em}\times \frac{\partial P(\mathbf Y_A=\bm 0_A,\mathbf Y_B\leq \mathbf w_B,\mathbf Y_C=\bm 0_C)}{\partial \mathbf w_B} \frac{\partial P(\mathbf Y_A=\bm 0_A,\mathbf Y_B=\bm 0_B,\mathbf Y_C\leq \mathbf w_C)}{\partial \mathbf w_C},\hspace{3em}\mathbf w_A \in L_A^2,\mathbf w_B \in L_B^2, \mathbf w_C \in L_C^2.
\end{cases}
\end{aligned}
\]	
Since the support of $\mathbf Y=(\mathbf Y_A,\mathbf Y_B,\mathbf Y_C)$  is $\Omega =\left( L_A^2 \times \bm{0}_{B}\times \bm{0}_{C} \right) \cup  \left( \bm{0}_{A}\times L_B^2\times  \bm{0}_{C} \right) \cup  \left( \bm{0}_{A}\times \bm{0}_{B}\times L_C^2  \right)$, then $\partial P(\mathbf Y_A=\bm 0_A,\mathbf Y_B=\bm 0_B,\mathbf Y_C \leq \mathbf w_C)/\partial \mathbf w_C \neq 0$ for $\mathbf w_C \in L_C^2$, $\partial P(\mathbf Y_A=\bm 0_A,\mathbf Y_B \leq \mathbf w_B,\mathbf Y_C=\bm 0_C)/\partial \mathbf w_B \neq 0$ for $\mathbf w_B \in L_B^2$ and $\partial P(\mathbf Y_A\leq \mathbf w_A,\mathbf Y_B=\bm 0_B,\mathbf Y_C=\bm 0_C)/\partial \mathbf w_A\neq 0$ for $\mathbf w_A\in L_A^2$. Thus, from $f(\mathbf w_A,\mathbf w_B,\mathbf w_C)$, we 
can obtain that the support of $\mathbf W=(\mathbf W_A,\mathbf W_B,\mathbf W_C)$ is $[0,\infty)^{|A|+|B|+|C|}$, which contains  $\Omega$.  Based on the definition of conditional density in \citet[\S 7.4]{shorack2017probability}, if  $\mathbf w_A 
\in L_A^2$, $\mathbf w_B \in L_B^2$ and  $\mathbf w_C \in L_C^2$, then we have that the conditional densities at $(\mathbf w_A,\bm 0_B,\bm 0_C)$, $(\bm 0_A,\mathbf w_B,\bm 0_C)$ and  $(\bm 0_A,\bm 0_B,\mathbf w_C)$,  given that all the points are on the curve of $\Omega$ are
\[
\begin{aligned}
&\frac{f(\mathbf w_A,\bm 0_B,\bm 0_C)}{\int_{L_A^2} f(\mathbf w_A,\bm 0_B,\bm 0_C) {\rm d} \mathbf w_A  + \int_{L_B^2} f(\bm 0_A, \mathbf w_B, \bm 0_C) {\rm d} \mathbf w_B + \int_{L_C^2} f(\bm 0_A, \bm 0_B, \mathbf w_C) {\rm d} \mathbf w_C}\\
=&\frac{p_A(1-p_B)(1-p_C)\frac{\partial P(\mathbf Y_A\leq \mathbf w_A, \mathbf Y_B=\bm 0_B, \mathbf Y_C=\bm 0_C )}{\partial \mathbf w_A}
}{\alpha_1\{p_A(1-p_B)(1-p_C)+ (1-p_A)p_B(1-p_C)+(1-p_A)(1-p_B)p_C \}   },
\end{aligned}
\]
\[
\begin{aligned}
&\frac{f(\bm 0_A,\mathbf w_B,\bm 0_C)}{ \int_{L_A^2} f(\mathbf w_A,\bm 0_B,\bm 0_C) {\rm d} \mathbf w_A  + \int_{L_B^2} f(\bm 0_A, \mathbf w_B,\bm 0_C) {\rm d} \mathbf w_B + \int_{L_C^2} f(\bm 0_A, \bm 0_B, \mathbf w_C) {\rm d} \mathbf w_C}\\
=&\frac{(1-p_A)p_B(1-p_C)\frac{\partial P(\mathbf Y_A=\bm 0_A,\mathbf Y_B\leq \mathbf w_B,\mathbf Y_C=\bm 0_C)}{\partial \mathbf w_B} }{\alpha_2\{p_A(1-p_B)(1-p_C)+ (1-p_A)p_B(1-p_C)+(1-p_A)(1-p_B)p_C \}   },
\end{aligned}
\]
and 
\[
\begin{aligned}
&\frac{f(\bm 0_A,\bm 0_B,\mathbf w_C)}{ \int_{L_A^2} f(\mathbf w_A,\bm 0_B,\bm 0_C) {\rm d} \mathbf w_A  + \int_{L_B^2} f(\bm 0_A, \mathbf w_B,\bm 0_C) {\rm d} \mathbf w_B+ \int_{L_C^2} f(\bm 0_A, \bm 0_B, \mathbf w_C) {\rm d} \mathbf w_C }\\
=&\frac{(1-p_A)(1-p_B)p_C\frac{\partial P(\mathbf Y_A=\bm 0_A,\mathbf Y_B=\bm 0_B,\mathbf Y_C\leq \mathbf w_C)}{\partial \mathbf w_C} }{(1-\alpha_1-\alpha_2)\{p_A(1-p_B)(1-p_C)+ (1-p_A)p_B(1-p_C)+(1-p_A)(1-p_B)p_C \}   }.
\end{aligned}
\]
Then we can get the value of $P(\mathbf W_A\leq \mathbf w_A, \mathbf W_B\leq \mathbf w_B, \mathbf W_C\leq \mathbf w_C\mid \mathbf W_A\in L_A^2,\mathbf W_B=\bm{0}_{B},\mathbf W_C=\bm{0}_{C} \,\text{or}\, \mathbf W_A=\bm{0}_{A}, \mathbf W_B\in L_B^2, \mathbf W_C=\bm{0}_{C}\,\text{or}\,\mathbf W_A=\bm{0}_{A}, \mathbf W_B=\bm{0}_{B}, \mathbf W_C\in L_C^2)=P(\mathbf W_A\leq \mathbf w_A, \mathbf W_B\leq \mathbf w_B, \mathbf W_C\leq \mathbf w_C\mid \mathbf W\in \Omega)$. When $\mathbf w_A \in L_A^1,\mathbf w_B \in L_B^1, \mathbf w_C \in L_C^1$,
\[
P(\mathbf W_A\leq \mathbf w_A, \mathbf W_B\leq \mathbf w_B, \mathbf W_C\leq \mathbf w_C\mid \mathbf W\in \Omega)
=0=P(\mathbf Y_A\leq \mathbf w_A, \mathbf Y_B\leq \mathbf w_B, \mathbf Y_C\leq \mathbf w_C).
\]
When $\mathbf w_A \in L_A^1,\mathbf w_B \in L_B^1, \mathbf w_C \in L_C^2$,	
\[
\begin{aligned}
  &P(\mathbf W_A\leq \mathbf w_A, \mathbf W_B\leq \mathbf w_B, \mathbf W_C\leq \mathbf w_C\mid \mathbf W\in \Omega)\\
  =&\frac{\int _{[\bm{0}_C,\mathbf w_C]\setminus[0,1]^{|C|} } f(\bm 0_A,\bm 0_B,\mathbf t){\rm d} \mathbf t}{ \int_{L_A^2} f(\mathbf w_A,\bm 0_B,\bm 0_C) {\rm d} \mathbf w_A  + \int_{L_B^2} f(\bm 0_A, \mathbf w_B,\bm 0_C) {\rm d} \mathbf w_B + \int_{L_C^2} f(\bm 0_A, \bm 0_B, \mathbf w_C) {\rm d} \mathbf w_C }\\
    =&\frac{(1-p_A)(1-p_B)p_C P(\mathbf Y_A=\bm 0_A,\mathbf Y_B=\bm 0_B,\mathbf Y_C\leq \mathbf w_C)}{(1-\alpha_1-\alpha_2)\{p_A(1-p_B)(1-p_C)+ (1-p_A)p_B(1-p_C)+(1-p_A)(1-p_B)p_C \}}\\
        =&\frac{(1-p_A)(1-p_B)p_C P(\mathbf Y_A\leq \mathbf w_A,\mathbf Y_B\leq \mathbf w_B,\mathbf Y_C\leq \mathbf w_C)}{(1-\alpha_1-\alpha_2)\{p_A(1-p_B)(1-p_C)+ (1-p_A)p_B(1-p_C)+(1-p_A)(1-p_B)p_C \}}.
\end{aligned}
\]
When $\mathbf w_A \in L_A^1,\mathbf w_B \in L_B^2, \mathbf w_C \in L_C^1$,	
\[
\begin{aligned}
  &P(\mathbf W_A\leq \mathbf w_A, \mathbf W_B\leq \mathbf w_B, \mathbf W_C\leq \mathbf w_C\mid \mathbf W\in \Omega)\\
  =&\frac{\int _{[\bm{0}_B,\mathbf w_B]\setminus[0,1]^{|B|} } f(\bm{0}_A,\mathbf t,\bm{0}_C){\rm d} \mathbf t}{ \int_{L_A^2} f(\mathbf w_A,\bm{0}_B,\bm{0}_C) {\rm d} \mathbf w_A  + \int_{L_B^2} f(\bm{0}_A, \mathbf w_B,\bm{0}_C) {\rm d} \mathbf w_B + \int_{L_C^2} f(\bm{0}_A, \bm{0}_B, \mathbf w_C) {\rm d} \mathbf w_C }\\
    =&\frac{(1-p_A)p_B(1-p_C)P(\mathbf Y_A= \bm{0}_A,\mathbf Y_B\leq \mathbf w_B,\mathbf Y_C=\bm{0}_C)}{\alpha_2\{p_A(1-p_B)(1-p_C)+ (1-p_A)p_B(1-p_C)+(1-p_A)(1-p_B)p_C \}}\\
        =&\frac{(1-p_A)p_B(1-p_C)P(\mathbf Y_A\leq \mathbf w_A,\mathbf Y_B\leq \mathbf w_B,\mathbf Y_C\leq \mathbf w_C)}{\alpha_2\{p_A(1-p_B)(1-p_C)+ (1-p_A)p_B(1-p_C)+(1-p_A)(1-p_B)p_C \}}.
\end{aligned}
\]
When $\mathbf w_A \in L_A^1,\mathbf w_B \in L_B^2, \mathbf w_C \in L_C^2$,	due to $P(\mathbf Y_A\leq \mathbf w_A,\mathbf Y_B\leq \mathbf w_B,\mathbf Y_C\leq \mathbf w_C)=P(\mathbf Y_A=\bm 0_A,\mathbf Y_B\leq \mathbf w_B,\mathbf Y_C=\bm 0_C)+P(\mathbf Y_A=\bm 0_A,\mathbf Y_B=\bm 0_B,\mathbf Y_C\leq \mathbf w_C)$,
\[
\begin{aligned}
  &P(\mathbf W_A\leq \mathbf w_A, \mathbf W_B\leq \mathbf w_B, \mathbf W_C\leq \mathbf w_C\mid \mathbf W\in \Omega)\\
  =&\frac{\int _{[\bm{0}_B,\mathbf w_B]\setminus[0,1]^{|B|} } f(\bm{0}_A,\mathbf t,\bm{0}_C){\rm d} \mathbf t+\int _{[\bm{0}_C,\mathbf w_C]\setminus[0,1]^{|C|} } f(\bm{0}_A,\bm{0}_B,\mathbf t){\rm d} \mathbf t}{  \int_{L_A^2} f(\mathbf w_A,\bm{0}_B,\bm{0}_C) {\rm d} \mathbf w_A  + \int_{L_B^2} f(\bm{0}_A, \mathbf w_B,\bm{0}_C) {\rm d} \mathbf w_B + \int_{L_C^2} f(\bm{0}_A, \bm{0}_B, \mathbf w_C) {\rm d} \mathbf w_C }\\
        =&\frac{\{(1-p_A)p_B(1-p_C)(1-\alpha_1-\alpha_2)-(1-p_A)(1-p_B)p_C\alpha_2\}P(\mathbf Y_A=\bm{0}_A,\mathbf Y_B\leq \mathbf w_B,\mathbf Y_C=\bm{0}_C)}{\alpha_2(1-\alpha_1-\alpha_2)\{p_A(1-p_B)(1-p_C)+ (1-p_A)p_B(1-p_C)+(1-p_A)(1-p_B)p_C \}}\\
       &+ \frac{(1-p_A)(1-p_B)p_C P(\mathbf Y_A\leq \mathbf w_A,\mathbf Y_B\leq \mathbf w_B,\mathbf Y_C\leq \mathbf w_C)}{(1-\alpha_1-\alpha_2)\{p_A(1-p_B)(1-p_C)+ (1-p_A)p_B(1-p_C)+(1-p_A)(1-p_B)p_C \}}.
\end{aligned}
\]
When $\mathbf w_A \in L_A^2,\mathbf w_B \in L_B^1, \mathbf w_C \in L_C^1$,	
\[
\begin{aligned}
  &P(\mathbf W_A\leq \mathbf w_A, \mathbf W_B\leq \mathbf w_B, \mathbf W_C\leq \mathbf w_C\mid \mathbf W\in \Omega)\\
  =&\frac{\int _{[\bm{0}_A,\mathbf w_A]\setminus[0,1]^{|A|} } f(\mathbf t,\bm{0}_B,\bm{0}_C){\rm d} \mathbf t}{ \int_{L_A^2} f(\mathbf w_A,\bm{0}_B,\bm{0}_C) {\rm d} \mathbf w_A  + \int_{L_B^2} f(\bm{0}_A, \mathbf w_B,\bm{0}_C) {\rm d} \mathbf w_B + \int_{L_C^2} f(\bm{0}_A, \bm{0}_B, \mathbf w_C) {\rm d} \mathbf w_C }\\
    =&\frac{p_A(1-p_B)(1-p_C)P(\mathbf Y_A\leq \mathbf w_A,\mathbf Y_B=\bm{0}_B,\mathbf Y_C=\bm{0}_C)}{\alpha_1\{p_A(1-p_B)(1-p_C)+ (1-p_A)p_B(1-p_C)+(1-p_A)(1-p_B)p_C \}}\\
        =&\frac{p_A(1-p_B)(1-p_C)P(\mathbf Y_A\leq \mathbf w_A,\mathbf Y_B\leq \mathbf w_B,\mathbf Y_C\leq \mathbf w_C)}{\alpha_1\{p_A(1-p_B)(1-p_C)+ (1-p_A)p_B(1-p_C)+(1-p_A)(1-p_B)p_C \}}.
\end{aligned}
\]
When $\mathbf w_A \in L_A^2,\mathbf w_B \in L_B^1, \mathbf w_C \in L_C^2$,	due to $P(\mathbf Y_A\leq \mathbf w_A,\mathbf Y_B\leq \mathbf w_B,\mathbf Y_C\leq \mathbf w_C)=P(\mathbf Y_A\leq \mathbf w_A,\mathbf Y_B=\bm 0_B,\mathbf Y_C=\bm 0_C)+P(\mathbf Y_A=\bm 0_A,\mathbf Y_B=\bm 0_B,\mathbf Y_C\leq \mathbf w_C)$,
\[
\begin{aligned}
  &P(\mathbf W_A\leq \mathbf w_A, \mathbf W_B\leq \mathbf w_B, \mathbf W_C\leq \mathbf w_C\mid \mathbf W\in \Omega)\\
  =&\frac{\int _{[\bm{0}_A,\mathbf w_A]\setminus[0,1]^{|A|} } f(\mathbf t,\bm{0}_B,\bm{0}_C){\rm d} \mathbf t+\int _{[\bm{0}_C,\mathbf w_C]\setminus[0,1]^{|C|} } f(\bm{0}_A,\bm{0}_B,\mathbf t){\rm d} \mathbf t}{  \int_{L_A^2} f(\mathbf w_A,\bm{0}_B,\bm{0}_C) {\rm d} \mathbf w_A  + \int_{L_B^2} f(\bm{0}_A, \mathbf w_B,\bm{0}_C) {\rm d} \mathbf w_B + \int_{L_C^2} f(\bm{0}_A, \bm{0}_B, \mathbf w_C) {\rm d} \mathbf w_C }\\
        =&\frac{\{p_A(1-p_B)(1-p_C)(1-\alpha_1-\alpha_2)-(1-p_A)(1-p_B)p_C\alpha_1\}P(\mathbf Y_A\leq \mathbf w_A,\mathbf Y_B=\bm{0}_B,\mathbf Y_C=\bm{0}_C)}{\alpha_1(1-\alpha_1-\alpha_2)\{p_A(1-p_B)(1-p_C)+ (1-p_A)p_B(1-p_C)+(1-p_A)(1-p_B)p_C \}}\\
       &+ \frac{(1-p_A)(1-p_B)p_C P(\mathbf Y_A\leq \mathbf w_A,\mathbf Y_B\leq \mathbf w_B,\mathbf Y_C\leq \mathbf w_C)}{(1-\alpha_1-\alpha_2)\{p_A(1-p_B)(1-p_C)+ (1-p_A)p_B(1-p_C)+(1-p_A)(1-p_B)p_C \}}.
\end{aligned}
\]
When $\mathbf w_A \in L_A^2,\mathbf w_B \in L_B^2, \mathbf w_C \in L_C^1$,	due to $P(\mathbf Y_A\leq \mathbf w_A,\mathbf Y_B\leq \mathbf w_B,\mathbf Y_C\leq \mathbf w_C)=P(\mathbf Y_A\leq \mathbf w_A,\mathbf Y_B=\bm 0_B,\mathbf Y_C=\bm 0_C)+P(\mathbf Y_A=\bm 0_A,\mathbf Y_B\leq \mathbf w_B,\mathbf Y_C=\bm 0_C)$,
\[
\begin{aligned}
  &P(\mathbf W_A\leq \mathbf w_A, \mathbf W_B\leq \mathbf w_B, \mathbf W_C\leq \mathbf w_C\mid \mathbf W\in \Omega)\\
  =&\frac{\int _{[\bm{0}_A,\mathbf w_A]\setminus[0,1]^{|A|} } f(\mathbf t,\bm{0}_B,\bm{0}_C){\rm d} \mathbf t+\int _{[\bm{0}_B,\mathbf w_B]\setminus[0,1]^{|B|} } f(\bm{0}_A,\mathbf t,\bm{0}_C){\rm d}\mathbf t}{  \int_{L_A^2} f(\mathbf w_A,\bm{0}_B,\bm{0}_C) {\rm d} \mathbf w_A  + \int_{L_B^2} f(\bm{0}_A, \mathbf w_B,\bm{0}_C) {\rm d} \mathbf w_B + \int_{L_C^2} f(\bm{0}_A, \bm{0}_B, \mathbf w_C) {\rm d} \mathbf w_C }\\
        =&\frac{\{p_A(1-p_B)(1-p_C)\alpha_2-(1-p_A)p_B(1-p_C)\alpha_1\}P(\mathbf Y_A\leq \mathbf w_A,\mathbf Y_B=\bm{0}_B,\mathbf Y_C=\bm{0}_C)}{\alpha_1\alpha_2\{p_A(1-p_B)(1-p_C)+ (1-p_A)p_B(1-p_C)+(1-p_A)(1-p_B)p_C \}}\\
       &+ \frac{(1-p_A)p_B(1-p_C)P(\mathbf Y_A\leq \mathbf w_A,\mathbf Y_B\leq \mathbf w_B,\mathbf Y_C\leq \mathbf w_C)}{\alpha_2\{p_A(1-p_B)(1-p_C)+ (1-p_A)p_B(1-p_C)+(1-p_A)(1-p_B)p_C \}}.
\end{aligned}
\]
When $\mathbf w_A \in L_A^2,\mathbf w_B \in L_B^2, \mathbf w_C \in L_C^2$,	due to $P(\mathbf Y_A\leq \mathbf w_A,\mathbf Y_B\leq \mathbf w_B,\mathbf Y_C\leq \mathbf w_C)=P(\mathbf Y_A\leq \mathbf w_A,\mathbf Y_B=\bm 0_B,\mathbf Y_C=\bm 0_C)+P(\mathbf Y_A=\bm 0_A,\mathbf Y_B\leq \mathbf w_B,\mathbf Y_C=\bm 0_C)+P(\mathbf Y_A=\bm 0_A,\mathbf Y_B=\bm 0_B,\mathbf Y_C\leq \mathbf w_C)$,
\begin{multline*}
  P(\mathbf W_A\leq \mathbf w_A, \mathbf W_B\leq \mathbf w_B, \mathbf W_C\leq \mathbf w_C\mid \mathbf W\in \Omega)\\
  =\frac{\int _{[\bm{0}_A,\mathbf w_A]\setminus[0,1]^{|A|} } f(\mathbf t,\bm{0}_B,\bm{0}_C){\rm d} \mathbf t+\int _{[\bm{0}_B,\mathbf w_B]\setminus[0,1]^{|B|} } f(\bm{0}_A,\mathbf t,\bm{0}_C){\rm d} \mathbf t+\int _{[\bm{0}_C,\mathbf w_C]\setminus[0,1]^{|C|} } f(\bm{0}_A,\bm{0}_B,\mathbf t){\rm d} \mathbf t}{  \int_{L_A^2} f(\mathbf w_A,\bm{0}_B,\bm{0}_C) {\rm d} \mathbf w_A  + \int_{L_B^2} f(\bm{0}_A, \mathbf w_B,\bm{0}_C) {\rm d} \mathbf w_B + \int_{L_C^2} f(\bm{0}_A, \bm{0}_B, \mathbf w_C) {\rm d} \mathbf w_C }\\
        =\Bigg[\frac{\{p_A(1-p_B)(1-p_C)(1-\alpha_1-\alpha_2)-(1-p_A)(1-p_B)p_C\alpha_1\}P(\mathbf Y_A\leq \mathbf w_A,\mathbf Y_B=\bm{0}_B,\mathbf Y_C=\bm{0}_C)}{\alpha_1(1-\alpha_1-\alpha_2)}\\
      + \frac{\{(1-p_A)p_B(1-p_C)(1-\alpha_1-\alpha_2)-(1-p_A)(1-p_B)p_C\alpha_2\}P(\mathbf Y_A=\bm{0}_A,\mathbf Y_B\leq \mathbf w_B,\mathbf Y_C=\bm{0}_C)}{\alpha_2(1-\alpha_1-\alpha_2)}\\
      +\frac{(1-p_A)(1-p_B)p_C P(\mathbf Y_A\leq \mathbf w_A,\mathbf Y_B\leq \mathbf w_B,\mathbf Y_C\leq \mathbf w_C)}{1-\alpha_1-\alpha_2}\Bigg]\\
      \times \frac{1}{p_A(1-p_B)(1-p_C)+ (1-p_A)p_B(1-p_C)+(1-p_A)(1-p_B)p_C}.
\end{multline*}
By comparison, we can find that $P(\mathbf W_A\leq \mathbf w_A, \mathbf W_B\leq \mathbf w_B, \mathbf W_C\leq \mathbf w_C\mid \mathbf W\in \Omega)=P(\mathbf Y_A\leq \mathbf w_A, \mathbf Y_B\leq \mathbf w_B, \mathbf Y_C\leq \mathbf w_C)$ if and only if 
\begin{equation}
\label{eqn:conditon2}
\begin{cases}
\alpha_1=\frac{p_A(1-p_B)(1-p_C)}{p_A(1-p_B)(1-p_C)+ (1-p_A)p_B(1-p_C)+(1-p_A)(1-p_B)p_C},\\
\alpha_2=\frac{(1-p_A)p_B(1-p_C)}{p_A(1-p_B)(1-p_C)+ (1-p_A)p_B(1-p_C)+(1-p_A)(1-p_B)p_C}.
\end{cases}
\end{equation}
We have hence shown that under condition \eqref{eqn:conditon2}, $(\mathbf W \mid \mathbf W\in \Omega) \,{\buildrel d \over =}\, \mathbf Y$ if the density function of $F(\mathbf w_A,\mathbf w_B,\mathbf w_C)$ exists.

If the density function of  $F(\mathbf w_A,\mathbf w_B,\mathbf w_C)$ does not exist, from the definition of conditional probability in \citet[\S 7.4]{shorack2017probability}, we know that the above formulas for $P(\mathbf W_A\leq \mathbf w_A, \mathbf W_B\leq \mathbf w_B, \mathbf W_C\leq \mathbf w_C\mid \mathbf W\in \Omega)$ still satisfy the condition for conditional probability. Thus, we still can define $P(\mathbf W_A\leq \mathbf w_A, \mathbf W_B\leq \mathbf w_B, \mathbf W_C\leq \mathbf w_C\mid \mathbf W\in \Omega)$ in the same way as before. Under condition \eqref{eqn:conditon2}, we also have $P(\mathbf W_A\leq \mathbf w_A, \mathbf W_B\leq \mathbf w_B, \mathbf W_C\leq \mathbf w_C\mid \mathbf W\in \Omega)=P(\mathbf Y_A\leq \mathbf w_A,\mathbf Y_B\leq \mathbf w_B,\mathbf Y_C\leq \mathbf w_C)$. Then $(\mathbf W \mid \mathbf W\in \Omega) \,{\buildrel d \over =}\, \mathbf Y$.

	In summary, according to the definition of $\mathbf Y_A \ind_o \mathbf Y_C \mid  \mathbf Y_B$, we have $\mathbf Y_A \ind_o \mathbf Y_C \mid  \mathbf Y_B$.
\end{proof}

	\thispagestyle{empty}
	\bibliographystyle{apalike}
	\bibliography{causal}

\end{document}